\documentclass[letterpaper]{IEEETran}

\usepackage[paper=letterpaper,top=0.60in,bottom=0.55in,right=0.47in,left=0.52in]{geometry}
\usepackage[comma,numbers,square,sort&compress]{natbib}
\usepackage{algorithm} %format of the algorithm
\usepackage{algorithmic} %format of the algorithm
\usepackage{multirow} %multirow for format of table
\usepackage{xcolor}
\usepackage{amsmath}
\usepackage{amsmath}
\usepackage{multirow}
\usepackage{latexsym}
\usepackage{setspace}
\usepackage{url}
\usepackage{subfigure}
\usepackage{fancyhdr}
\usepackage{amsmath}
\usepackage{multirow}
\usepackage{amssymb }
\usepackage{color}
\usepackage{graphics} % for pdf,  bitmapped graphics files
\usepackage{graphicx}
\usepackage{pmat}
\usepackage{arydshln}

\newtheorem{theorem}{\textbf{Theorem}}
\newtheorem{lemma}{\textbf{Lemma}}
\newtheorem{corollary}{\textbf{Corollary}}
\newtheorem{remark}{\textbf{Remark}}
\newtheorem{definition}{\textbf{Definition}}
\newtheorem{proposition}{\textbf{Proposition}}

\newtheorem{assumption}{\textbf{Assumption}}
\newtheorem{problem}{\textbf{Problem}}

\begin{document}
\fontsize{10pt}{11pt}\selectfont

\title{\LARGE \bf Structural Controllability of a Networked Dynamic System with LFT Parameterized Subsystems}

\author{Yuan Zhang and Tong Zhou$^{\dag}$
\thanks{This work was supported in part by the NNSFC under Grant  61733008 and 61573209. {\emph{(Corresponding author: Tong Zhou.)$^{\dag}$}}}
\thanks{Yuan Zhang and Tong Zhou are with the Department of Automation, Tsinghua University, Beijing, 100084, P. R. China.
        {(Email: {\tt\small zhangyuan14@mails.tsinghua.edu.cn, tzhou@mail.tsinghua.edu.cn}.)}}
}

\maketitle

\begin{abstract}
This paper studies structural controllability for a networked dynamic system (NDS), in which each subsystem may have different dynamics, and unknown parameters may exist both in subsystem dynamics and in subsystem interconnections. In addition, subsystem parameters are parameterized by a linear fractional transformation (LFT). It is proven that controllability keeps to be a generic property for this kind of NDSs. Some necessary and sufficient conditions are then established respectively for them to be structurally controllable, to have a fixed uncontrollable mode, and to have a parameter dependent uncontrollable mode, under the condition that each subsystem interconnection link can take a weight independently. These conditions are scalable, and in their verifications, all arithmetic calculations are performed separately on each subsystem. In addition, these conditions also reveal influences on NDS controllability from subsystem input-output relations, subsystem uncontrollable modes and subsystem interconnection topology. Based on these observations, the problem of selecting the minimal number of subsystem interconnection links is studied under the requirement of constructing a structurally controllable NDS. A heuristic method is derived with some provable approximation bounds and a low computational complexity.
\end{abstract}  %the weights of the subsystem interconnection links can take values freely %nonzero entries of the subsystem connection matrix can take values freely % under the condition that all subsystem first principle parameters and weights of subsystem interconnection links can take values freely


\section{Introduction}

Controllability and observability are now relatively classic concepts in system analysis and synthesis, which are closely related to stabilization,  existence of an optimal control, convergence of state estimators and some other fundamental issues  \cite{Zhou_robust_book}, \cite{Kailath_1980}, \cite{Modern_Control_Ogata}. The past two decades have seen a renewed  research interest in controllability and observability  with the emergence of complex networks, such as biological transduction networks \cite{nature}, power networks \cite{FBullo2013} and social networks \cite{nature}, \cite{controllability metrics}.
Recent applications of controllability and observability include controlling a biological network using as less inputs as possible \cite{nature}, attack detection in cyber-physical systems \cite{FBullo2013}, achieving consensus in multi-agent systems \cite{partitions},  distributed control of power generation networks \cite{controllability metrics}, etc.

Apart from the extensively adopted concept of controllability proposed by Kalman \cite{Kailath_1980}, structural controllability proposed by Lin \cite{Lin_1974} is also widely utilized in which each entry of system matrices is either fixed to be zero or allowed to take an arbitrary real value.
A system with this kind of system matrices is called a structured system, and is said to be  weakly structurally controllable, which is sometimes called structurally controllable, if there exists one numerical realization such that the associated system is controllable. A closely related concept is strong structural controllability \cite{strong_controllability}, which requires that each numerical realization corresponds to a controllable system. While these two concepts significantly differ from each other in their engineering motivations, there do not exist great differences in applying them to an actual engineering problem, noting that controllability is a generic system property, which means that, structural controllability of a structured system guarantees that all its numerical realizations,  except for a set of zero Lebesgue measure in the associated parameter space, correspond to a controllable system \cite{generic}. A prominent characteristic of structural controllability is that its criteria often have a clear graphical interpretation, which explicitly reveals information flows in a dynamic system, and can usually be easily verified graphically \cite{generic}, \cite{nature}, \cite{Composability}. {{Recently, structural controllability has also been extended to some nonlinear systems, for example, the driftless bilinear systems \cite{Tsopelakos bilinear}}.}

When controllability of a networked dynamic system (NDS) is to be investigated,  some challenging issues arise due to nodal dynamics and high dimensions of its global system matrices \cite{plos one}, \cite{zhou_2015}. Recall that when a criterion for a lumped plant is applied to an NDS with a high dimensional state vector, numerical instability and/or computational prohibitiveness may emerge. Owing to efforts from many researchers, various graphical criteria for (structural) controllability have been  obtained \cite{nature}, \cite{Composability}, \cite{partitions}, \cite{generic}. These criteria, mainly focusing on an NDS with each of its subsystems being modelled as a first order differential equation, give many intuitive and useful insights on how the topology influences performances of an NDS. For example, it is shown in \cite{nature, plos one} that the minimum number of driver nodes needed to make a network controllable is closely related to its degree distribution, in case that the number of self-loops is negligible compared to the network size. Recently, extensive attentions have been moved to controllability of a compositional networked system in which each subsystem may have high order dynamics. For example, controllability is discussed for a network-of-network system via Cartesian product in \cite{network of networks}. In \cite{zhou_2015, zhou_2018, Y_Zhang_2016}, a general NDS model is adopted in which each subsystem can have different dynamics and the subsystem interconnection is arbitrary.

On the other hand, various results have been reported in which more general parameter interdependencies are adopted than those of the structured system model. More precisely, a linear parametrization is adopted in \cite{Morse_1976} which assumes that each entry of system matrices is an affine function of some free parameters.  A so called ``matrix net" which is a generalization of matrix pencil, is used in  \cite{Anderson_1982}. In \cite{Murota_Book}, the author introduces a concept called ``mixed matrix", where the nonzero entry of a matrix could be either a fixed constant or a free parameter, and uses matroid to study structural controllability under some assumptions. Recent related works include graphical interpretations of the conditions of \cite{Morse_1976} in \cite{Liu_Fei}, and leader selection for structured descriptor systems in \cite{Clark_input_selection}. The obtained results, which are often expressed through algebraic operations, however, are in general not computationally efficient for a large-scale NDS. {{Besides, parameter interdependencies adopted in these investigations are  not general enough to describe all actual plants, noting that entries in system matrices are usually rational functions of the parameters that govern its movements \cite{generic}, \cite{Modern_Control_Ogata}, which include parameters directly describing a physical process, a chemical process or a biological process, such as mass, temperature, concentration of a chemical element, etc.  All these parameters are called a first principle parameter in this paper for brevity.  An appropriate framework for describing matrices with rational function entries appears to be the  linear fractional transformation (LFT) widely adopted in robust control theory \cite{Zhou_robust_book}. }}

In this paper, we investigate structural controllability of an NDS in which each subsystem can have high-order and heterogeneous dynamics, and unknown parameters are allowed to exist in both subsystem dynamics and in subsystem interconnections. {{
We adopt an LFT to model subsystem parameter interdependencies, which enables describing a large class of plants whose unknown entries in system matrices are rational functions of the plant first principle parameters, and contains many other descriptions, such as the linear parametrization adopted in \cite{Morse_1976} and \cite{Willems_1986} as a special case. At first, necessary and sufficient conditions are established for structural controllability of an LFT parameterized plant under a diagonalization assumption. These results are obtained by characterizing the existences of a parameter dependent uncontrollable mode and a fixed uncontrollable mode, based on structure analysis of the associated transfer function matrices (TFMs) and the matroid theory. Compared to the results reported in \cite{Morse_1976}, \cite{Anderson_1982}, \cite{Willems_1986} and \cite{Murota_SIAM}, the adopted model is more general in describing relations between system matrices of a plant and its first principle parameters, and it has been made clear that the obtained conditions can be verified in polynomial time.}} Under the condition that each subsystem interconnection link can take a weight independently, some necessary and sufficient conditions are then derived for the NDS to be structurally controllable. These conditions can be verified efficiently, and the associated arithmetic operations are within each individual subsystem while the associated graphical operations are on the network topology, which makes them scalable and therefore attractive for large-scale NDSs. Moreover, these conditions explicitly illustrate how subsystem dynamics and subsystem interconnection topology, jointly influence controllability of an NDS. Furthermore, based on them, we consider the problem of designing a subsystem interaction topology that minimizes the number of interconnection links under structural controllability restriction. This problem is shown to be NP-hard, and a two-stage algorithm is suggested to approximate it with some provable sub-optimality guarantees. Additionally, we discuss the general computational complexity and hardness of structural controllability verification problems with more complicated parameter interdependencies.  % {{Additionally, we discuss the general computational complexity and hardness of structural controllability verification problems with more complicated parameter interdependencies, which indicates that a randomized algorithm may be feasible for the high-rank cases.}} %Each stage is accompanied with a sub-optimal solution guarantee, and the whole algorithm has some provable approximation bounds.  For $f(\lambda)\in F(\lambda)^{n_1\times n_2}$, $f(\lambda)\ne 0$ means that $f(\lambda)$ is not identically zero.

The rest of this paper is organized as follows. Section II provides problem statements and some preliminary results. Genericity is established in Section III for the controllability of an NDS with LFT parameterized subsystems. Necessary and sufficient conditions for structural controllability of LFT parameterized plants are presented in Section IV. These conditions are extended to NDSs in Section V taking their parameter structure into account. Section VI investigates designs of a minimal subsystem interaction topology to guarantee structural controllability for an NDS, while an illustrative example are given in Section VII. Finally,  Section VIII concludes this paper. Two appendices are included to provide some technical details.

\emph{Notations.}  Given two matrices $M = {\scriptsize { \begin{pmat}[{|}]
{{M_{11}}}& {{M_{12}}}\cr\-
{{M_{21}}}& {{M_{22}}}\cr
\end{pmat}}}$ and $ P$ with compatible dimensions,  if $I - {M_{22}} P$ is invertible, a lower LFT is defined as
${F_l}(M, P ) = {M_{11}} + {M_{12}} P {(I - {M_{22}} P )^{ - 1}}{M_{21}}$.  By ${\bf{{ diag}}}\{X_i|_{i=1}^N\}$ we denote the block diagonal matrix with its $i$-th diagonal block  being $X_i$, while ${\bf{ {{col}}}}\{X_i|_{i=1}^N\}$ the matrix stacked by $X_i|_{i=1}^N$ with its $i$-th row block being $X_i$. Denote by $F(\lambda)$ the field of all rational functions of $\lambda$ with real coefficients, and $F(\lambda)^{n_1\times n_2}$ the set of all $n_1\times n_2$ matrices with every entry in $F(\lambda)$. $\sigma(A)$ denotes the set of eigenvalues of a square matrix $A$.  $\mathbb{R}$, $\mathbb{C}$, $\mathbb{Z}$ and $\mathbb{N}$ denote the set of real, complex, integral and non-negative integral numbers, respectively.  Given $n\in \mathbb{N}$, define $[n]=\{1,2,...,n\}$.  For an $n_1\times n_2$ matrix $M$, $M_{ij}$ or $[M]_{ij}$ denotes its $(i,j)$-th entry, and for $J_1\subseteq [n_1]$ and $J_2\subseteq [n_2]$, $M_{J_1,J_2}$ denotes the submatrix of $M$ formed by rows indexed by $J_1$ and columns indexed by $J_2$, while $M_{J_2}$ the submatrix of $M$ formed by columns indexed by $J_2$. By $\{0,*\}^{n_1\times n_2}$ we denote an  $n_1\times n_2$ structured matrix, where $*$ denotes the entries which can take real values independently. By $||M||_0$ we denote the number of nonzero entries in a matrix $M$. Throughout this paper, a set of variables $S$ is said to be algebraically independent, if $S$ do not satisfy any non-trivial polynomial equation with its coefficients in $\mathbb{R}$.

\section{Problem Statement and Preliminaries}

\subsection{A Model for NDSs and Problem Statement}

In an actual NDS, it is often the case that its subsystems have high orders and heterogeneous dynamics. To describe dynamics of these NDSs, it appears convenient to utilize the spatially interconnected system model adopted in \cite{langbort_2004}, \cite{zhou_2015,zhou_2018, Y_Zhang_2016}, in which an NDS ${\bf{{\bf \Sigma}}}$ is constituted of $N$ subsystems, and the dynamics of its $i$-th subsystem ${\bf{\Sigma}}_i$ is described by
\begin{equation} \label{SubsystemDynamics}
\begin{array}{l}
\hspace*{-0.4cm}\left[\! \begin{matrix}
{{{\dot x}_i}(t)}\\
{{z_i}(t)}\\
{{y_i}(t)}
\end{matrix} \!\right] \!=\! \left\{\! {\left[\! \begin{matrix}
{A_{xx0}^{(i)}}&{A_{xv0}^{(i)}}&{B_{xu0}^{(i)}}\\
{A_{zx0}^{(i)}}&{A_{zv0}^{(i)}}&{B_{zu0}^{(i)}}\\
{C_{yx0}^{(i)}}&{C_{yv0}^{(i)}}&{D_{yu0}^{(i)}}
\end{matrix} \!\right] \!+ \! \left[\! \begin{matrix}
{E_1^{(i)}}\\
{E_2^{(i)}}\\
{E_3^{(i)}}
\end{matrix} \!\right]{ P ^{(i)}} \!\times \!} \right.\\
\hspace*{0.2cm}{{{(I - {H^{(i)}}{ P ^{(i)}})}^{ -1}}\left[\! \begin{matrix}
{F_1^{(i)}}&{F_2^{(i)}}&{F_3^{(i)}}
\end{matrix}\! \right]} {\Bigg{\}}}\left[ \begin{matrix} %\left.
{{x_i}(t)}\\
{{v_i}(t)}\\
{{u_i}(t)}
\end{matrix} \right]
\end{array}\end{equation}
where $t$ represents the temporal variable, {{$x_i(t)\in \mathbb{R}^{m_{xi}}$ is the state vector, $u_i(t)\in \mathbb{R}^{m_{ui}}$ is the external input vector, $y_i(t)\in \mathbb{R}^{m_{yi}}$ is the external output vector}}, $v_i(t)$ and $z_i(t)$ are respectively the signal received from other subsystems and the signal sent to other subsystems,  which are called the internal input and output vectors respectively, {{and whose dimensions can be greater than one}}.

In the above model, the system matrices of the $i$-th subsystem ${\bf{\Sigma}}_i$ are parameterized by a matrix $ P^{(i)}$, which is constituted from the unknown first principle parameters of this subsystem.  Moreover, $E_j^{(i)}$, $F_j^{(i)}$ and $H^{(i)}$, $j=1,2,3$, are adopted to represent some known constant matrices that reflect how a first principle parameter influences the subsystem matrices. The matrices with a subscript ``0'' represent the ingredient of the system matrices that does not vary with its first principle parameters. {{This forms an LFT parametrization widely adopted in robust control theory \cite{Zhou_robust_book}, and is able to describe system matrices which depend on its first principle parameters in a rational function way. }}
%A well encountered example of LFT parameterized plants is the mass/spring/damper system illustrated in Fig. 10.1 of \cite{Zhou_robust_book}.

\begin{figure} \centering
\subfigure[] { \label{fig:a}
\includegraphics[width=2.3in]{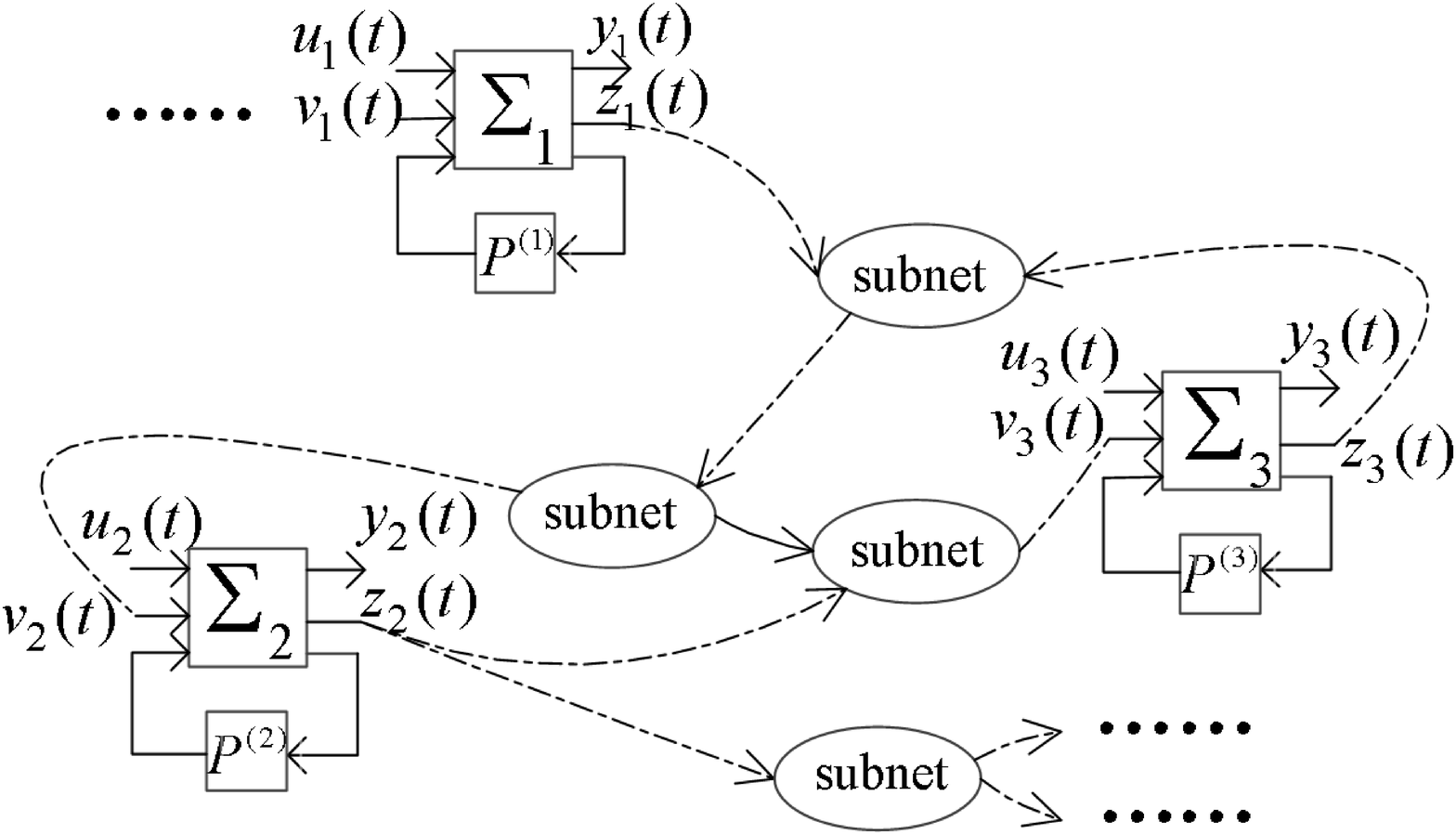}}
\subfigure[] { \label{fig:b}
\includegraphics[width=1.5in]{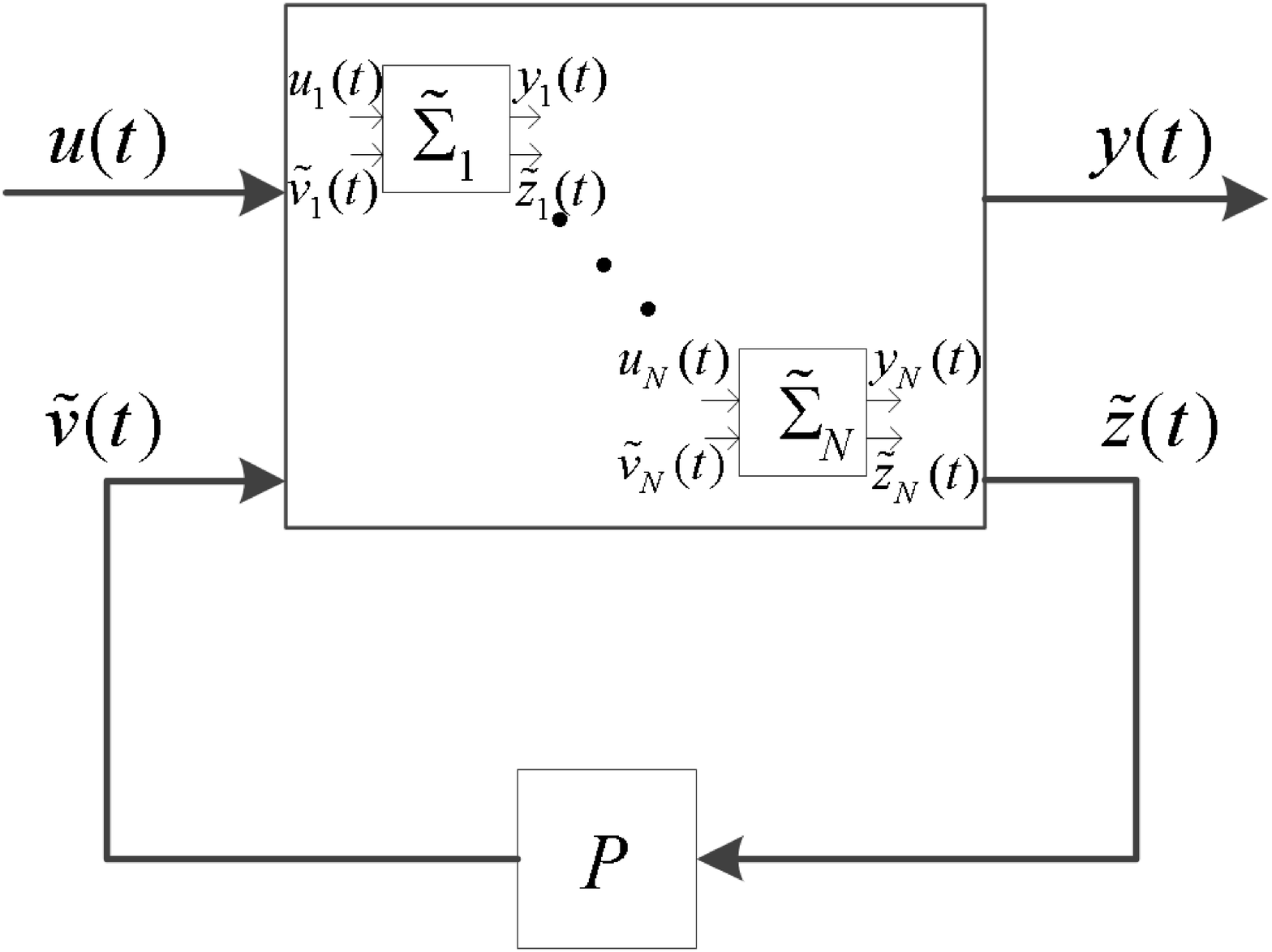}
}
\caption{\label{fig_network} (a): an illustrative example of the NDS ${\bf{\Sigma}}$; (b): an equivalent LFT scheme of the NDS ${\bf{\Sigma}}$.}
\end{figure}

In addition, interactions among subsystems are described by
\begin{equation}\label{SubsystemInteraction}
v(t)=\Phi z(t),
\end{equation}
where $v(t)={{\rm{{\bf{{col}}}}}}\{v_i(t)|_{i=1}^N\}$, $z(t)={{\rm{{\bf{{col}}}}}\{z_i(t)|_{i=1}^N\}}$. We call $\Phi$ the subsystem connection matrix (SCM), which describes the interconnection topology among subsystems. In actual networks, due to communication noises, {{inaccuracies of parameters describing the interaction channels}} or variations of spatial distances among subsystems, etc., the weights of the interconnection links may sometimes be hard to know exactly. Hence, while $\Phi$ is time-invariant, it is assumed that only the zero-nonzero patterns of $\Phi$ are known, that is, the positions of its elements that are fixed to be zero and those that are not constantly equal to zero. A fixed zero element in the SCM $\Phi$ means that the associated internal output of a subsystem does not directly affect another associated subsystem, while an element whose value is not constantly equal to zero means the contrary. In other words, the SCM $\Phi$ actually reflects the geometric structure of an NDS. With a little abuse of terminology, an element of the SCM $\Phi$ not fixed to be zero is also called a first principle parameter of the NDS in this paper. An NDS with the above dynamics is illustrated by Fig. \ref{fig_network}(a).

\begin{figure}
  \centering
  \includegraphics[width=3.5in]{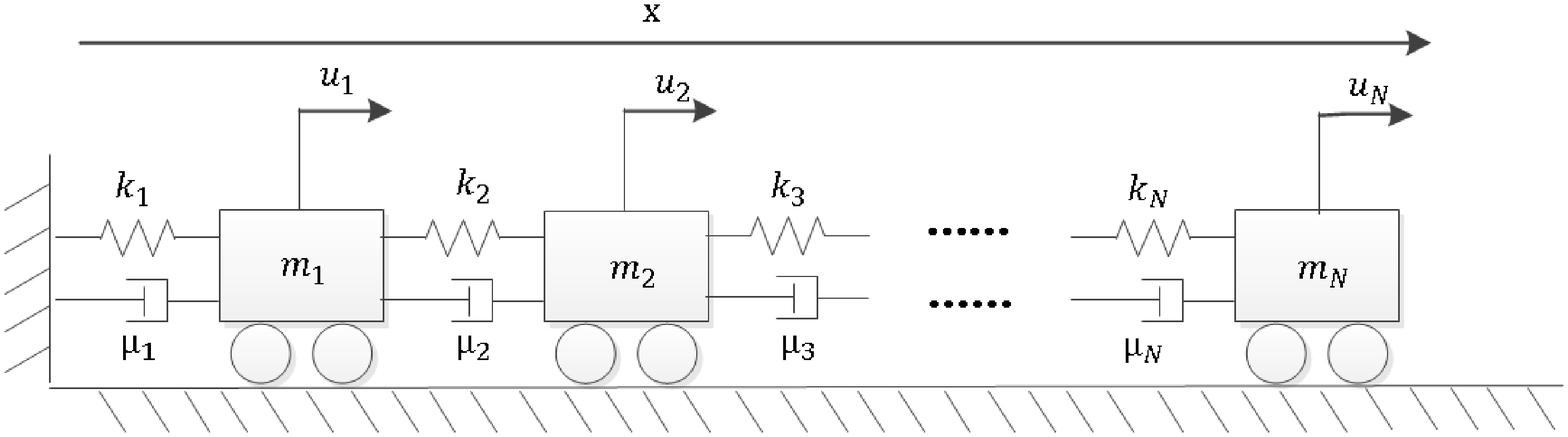}
  \caption{A vehicle-damper-spring chain system.}\label{damping_vehicle} %is composed by
\end{figure}

{{To illustrate the application significance of the NDS model, consider a mechanical system described in Fig. \ref{damping_vehicle} which consists of $N$ subsystems constituted from a vehicle, a damper and a spring. Obviously, the dynamics of each subsystem is determined by the mass of the vehicle, the constants of the springs and the dampers connected to it, in a rational function manner. Direct algebraic operations show that the system matrices of each subsystem in this plant can have the LFT form described by (\ref{SubsystemDynamics}), {{and dynamics of the whole system can be written as (\ref{SubsystemDynamics})-(\ref{SubsystemInteraction})}}. The details are omitted due to their straightforwardness.}} %Note that in this example, although some parameters may appear more than once in the associated SCM, they all have rank-one coefficient matrices.

The following assumption is adopted throughout this paper.

\begin{assumption}
The NDS ${\bf{{\bf \Sigma}}}$, as well as each of its subsystem ${\bf{{\bf \Sigma}}}_{i}$, $i\in [N]$, are well-posed for almost all feasible values of its first principle parameters.
\end{assumption}

The above assumption means that for each external input series ${{\rm{{\bf{{col}}}}}\{u_i(t)|_{i=1}^N\}}$, both the system states ${{\rm{{\bf{{col}}}}}\{x_i(t)|_{i=1}^N\}}$ and the external outputs ${{\rm{{\bf{{col}}}}}\{y_i(t)|_{i=1}^N\}}$ are uniquely determined for almost each $\Phi$ and each $P^{(i)}$, $i\in [N]$, having the corresponding prescribed structures. Under this assumption, the main problem discussed in this paper is stated as follows.

\begin{problem} {{Assume that all the nonzero entries of $P^{(1)},...,P^{(N)}$ and $\Phi$ are algebraically independent and time invariant}}, and except $P^{(i)}$, all the other values of the system matrices are prescribed for each subsystem ${\bf{{\bf \Sigma}}}_{i}$ of the NDS ${\bf{{\bf \Sigma}}}$, $i\in [N]$. Verify whether or not the NDS ${\bf{{\bf \Sigma}}}$ is structurally controllable. That is, whether or not there exist at least one feasible value respectively for $\Phi$ and each $ P^{(i)}$ for $i\in [N]$ with their prescribed structures, such that the corresponding NDS (\ref{SubsystemDynamics})-(\ref{SubsystemInteraction}) is controllable.
\end{problem}

In the next section, it is shown that if the NDS (\ref{SubsystemDynamics})-(\ref{SubsystemInteraction}) is structurally controllable, then for almost all feasible realizations of $ P^{(i)}|_{i=1}^N$ and $\Phi$, the corresponding NDS is controllable.

\begin{remark}
As mentioned earlier, LFT parametrization is capable of describing almost all rational function matrices \cite{Zhou_robust_book}.  It is worthwhile to mention that we restrict our attentions in this paper to the situation in which each first principle parameter appears only once in $P^{(i)}$, which can be directly extended to the case in which each first principle parameter has a rank-one coefficient matrix. This adoption is due to the following considerations. First, as discussed in Appendix A, the rank-one case is currently the most possible case for structural controllability verification that one can find deterministic algorithms with a polynomial time complexity or a sub-exponential time complexity. Second, a large class of traditional actual plants, which may be regarded as a subsystem of an NDS, satisfy this rank-one setting. These plants include many mechanical systems, electrical systems, as well as fluid systems \cite{Modern_Control_Ogata}, \cite{Zhou_robust_book}. The assumption that all nonzero entries of $\Phi$ are algebraically independent might be helpful, under some situations, in understanding the role of subsystem interconnection topology in controllability of NDSs.  %More specifically, many mechanical systems, electrical systems, as well as fluid systems can be modelled by an LFT parameterized state-space representation, with each of its first principle parameter appears only once in the diagonal of $P^{(i)}$ \cite{Modern_Control_Ogata}.  %More specifically, many mechanical systems, electrical systems, as well as fluid systems can be modelled by an LFT parameterized state-space representation, with each of its first principle parameter appears only once in the diagonal of $P^{(i)}$ \cite{Modern_Control_Ogata}.
%means that the weights of all subsystem interconnection links %most of the  is helpful to focusing on
\end{remark}

%\footnote{A transformation similar to (\ref{trans2LFT}) in Section V can transform the rank-one case to the case where each first principle parameter appears only once.}
\begin{remark} \label{remark3}
{{Various parametrizations have been proposed to describe system parameter interdependencies, e.g.,  the linear parametrization \cite{Morse_1976}, \cite{Willems_1986}, the matrix net \cite{Anderson_1982} and the mixed matrix descriptions \cite{Murota_Book}. It can be directly validated that these descriptions are special cases of the LFT parametrization ${F_l}([M_{ij}]_{i,j=1,2}, P)$ by setting $M_{22} \equiv 0$ and specifying the frequencies that the free parameters in $P$ appear. It is, however, worthwhile to emphasize that a nonzero $M_{22}$ is crucial and usually inevitable for enabling representation of rational functions by an LFT parametrization. }} %for the capacity of the LFT parametrizations in describing rational functions. %in  ${F_l}([M_{ij}]_{i,j=1,2}, P)$
\end{remark}

\vspace{0.1cm}
{{ \subsection{An Equivalent LFT Representation of the NDS}}}
To investigate structural controllability of the NDS described by (\ref{SubsystemDynamics})-(\ref{SubsystemInteraction}), its dynamics is rewritten as another LFT in which all first principle parameters in the model are included in one matrix.

For this purpose, two auxiliary internal input and output vectors are constructed for the $i$-th subsystem of the NDS (\ref{SubsystemDynamics})-(\ref{SubsystemInteraction}), which are denoted respectively by $v_i^{+}(t)$ and $z_i^{+}(t)$, and defined as
{\small{
\begin{equation} \label{auxiliary variable} \begin{split}
z_i^ + (t) &= [{F_1^{(i)}}\  {F_2^{(i)}}\  {F_3^{(i)}}]{\bf{{col}}}\{{{x_i}(t)},{{v_i}(t)}, {{u_i}(t)}\} + {H^{(i)}}v_i^ + (t),  \\
v_i^ + (t) &= { P ^{(i)}}z_i^ + (t).
\end{split}
\end{equation}}}Let ${\tilde z_i}(t)={{\rm{{\bf{{col}}}}}}\{z_i(t),z_i^{+}(t)\}$, ${\tilde v_i}(t) ={{\rm{{\bf{{col}}}}}}\{v_i(t),v_i^{+}(t)\}$. {{Denote the dimensions of $\tilde v_i(t)$ and $\tilde z_i(t)$ respectively by $m_{vi}$ and $m_{zi}$.}} Then it can be straightforwardly shown that, when the matrix $I-H^{(i)} P^{(i)}$ is invertible, the dynamics of Subsystem ${\bf{\Sigma}}_i$ can be equivalently expressed by (\ref{auxiliary variable}) and the following equation{\small{
\begin{equation}
\label{lft_sub}  \left[ \begin{matrix}
{{{\dot x}_i}(t)}\\
{{\tilde z_i}(t)}\\
{{y_i}(t)}
\end{matrix} \right] = \left[ \begin{matrix}
{A_{xx}^{(i)}}&{A_{xv}^{(i)}}&{B_{xu}^{(i)}}\\
{A_{zx}^{(i)}}&{A_{zv}^{(i)}}&{B_{zu}^{(i)}}\\
{C_{yx}^{(i)}}&{C_{yv}^{(i)}}&{D_{yu}^{(i)}}
\end{matrix} \right]\left[ \begin{matrix}
{x_{i}(t)}\\
{\tilde v_{i}(t)}\\
{u_{i}(t)}
\end{matrix} \right]
\end{equation}}}
in which{\small{
\begin{eqnarray*}
& & \hspace*{-0.8cm} A_{xx}^{(i)}=A_{xx0}^{(i)},\hspace{0.25cm}  A_{xv}^{(i)} = \left[ \begin{matrix}
{A_{xv0}^{(i)}}&{E_1^{(i)}}
\end{matrix} \right],\hspace{0.25cm}  A_{zx}^{(i)}{\rm{ = }}\left[ \begin{matrix}
{A_{zx0}^{(i)}}  \\
{F_1^{(i)}}
\end{matrix} \right], \\
& & \hspace*{-0.8cm}  A_{zv}^{(i)} = \left[ \begin{matrix}
{A_{zv0}^{(i)}}&{E_2^{(i)}}\\
{F_2^{(i)}}&{{H^{(i)}}}
\end{matrix} \right], \hspace{0.25cm} B_{xu}^{(i)}=B_{xu0}^{(i)}, \hspace{0.25cm} B_{zu}^{(i)} = \left[ \begin{matrix}
{B_{zu0}^{(i)}}\\
{F_3^{(i)}}
\end{matrix} \right], \\
& & \hspace*{-0.8cm} C_{yx}^{(i)} = C_{yx0}^{(i)}, \hspace{0.25cm} C_{yv}^{(i)} = \left[\begin{matrix}
{C_{yv0}^{(i)}}&{E_3^{(i)}}
\end{matrix} \right], \hspace{0.25cm}D_{yu}^{(i)} = D_{yu0}^{(i)}.
\end{eqnarray*}}}Partition the SCM $\Phi$ according to the dimensions of the vectors $v_i(t)$ and $z_j(t)$, $1\le i, j \le N$, and let $\Phi_{i,j}$ denote the $(i,j)$-th block of $\Phi$. Moreover, define vectors $\tilde v(t)$ and $\tilde z(t)$ respectively as $\tilde v(t) = {{\rm{{\bf{{col}}}}}}\{ {\tilde v_i}(t)|_{i=1}^N\} $ and $\tilde z(t) = {{\rm{{\bf{{col}}}}}}\{ {\tilde z_i}(t)|_{i=1}^N\}$.
Then the relation between $\tilde v(t)$ and $\tilde z(t)$ can be expressed as
\begin{equation}
\label{new_scl} \tilde v(t) =  P \tilde z(t),
\end{equation}
in which{\footnotesize{
\begin{equation}
\label{scm} P  = \begin{pmat}[{.|.||}]
{{\Phi _{1,1}}}&{}&{{\Phi _{1,2}}}&{}&{\cdots}&{{\Phi _{1,N}}}&{}\cr
{}&{ P _{}^{(1)}}&{}&0& \cdots &{}&0\cr\-
 \vdots & \vdots & \vdots & \vdots & \ddots & \vdots & \vdots \cr\-
{{\Phi _{N,1}}}&{}&{{\Phi _{N,2}}}&{}&{\cdots}&{{\Phi _{N,N}}}&{}\cr
{}&0&{}&0&{\cdots}&{}&{ P _{}^{(N)}}\cr
\end{pmat}.
\end{equation}}}

Equations (\ref{lft_sub})-(\ref{new_scl}) give an equivalent description for the input-output relations of the NDS ${\bf\Sigma}$, which
is illustrated by Fig. \ref{fig_network}(b). In the following, the system described by (\ref{lft_sub}) is called the augmented system of Subsystem ${\bf{\Sigma}}_i$, and is denoted by $\tilde{\bf{\Sigma}}_i$. {{It can be seen that the NDS (\ref{lft_sub})-(\ref{new_scl}) has the same subsystem interaction structure as that of the original NDS (\ref{SubsystemDynamics})-(\ref{SubsystemInteraction}),  except that some self-loops are introduced in (\ref{new_scl}).}}

Compared with (1)-(2) of \cite{zhou_2015}, it is obvious that the system described by (\ref{lft_sub})-(\ref{new_scl}) has completely the same form as that discussed in \cite{zhou_2015}. This implies that results of \cite{zhou_2015} can be directly applied to controllability analysis of an NDS with LFT parameterized subsystems. For this application, define matrices $A_{{\rm{\star }}*}$ and $B_{{\rm{\star }}*}$ with $\star, * = x, u, z$ or $v$ respectively as ${A_{{\rm{\star }} * }}{\rm{ = }} { \rm{{\bf{{\bf diag}}}}}\{ A_{{\rm{\star }} * }^{(i)}|_{i = 1}^N\}$, ${B_{{\rm{\star }} * }}{\rm{ = }}{\rm{{\bf{{\bf diag}}}}}\{B_{{\rm{\star }} * }^{(i)}|_{i = 1}^N\}$.  Moreover, define $M_{\#}=\sum\nolimits_{i=1}^N {m_{\#i}}$, with $\#=x,y,u,v$ and $z$. Then, a direct application of Theorem 1 in \cite{zhou_2015} leads to the following result, which gives a necessary and sufficient condition for the controllability of the NDS described by (\ref{SubsystemDynamics})-(\ref{SubsystemInteraction}) with a fixed $P$.

\begin{lemma}[\cite{zhou_2015}] \label{lemma_4}
Given a fixed $P$, assume that the associated NDS $\bf\Sigma$ (\ref{SubsystemDynamics})-(\ref{SubsystemInteraction}), as well as each of its subsystems,  is well-posed. Then this system is controllable, if and only if the following matrix valued polynomial (MVP) $M(\lambda)$ is of full row rank (FRR) for each $\lambda \in \mathbb{C}$,
\begin{equation}\label{second_MVP} M(\lambda ) = \left[\begin{matrix}
{\lambda I_{M_x} - {A_{{xx}}}}&{{B_{{xu}}}}&{ - {A_{xv}}{P}}\\
{ - {A_{zx}}}&{{B_{zu}}}&{I_{M_z} - {A_{zv}}{ P}}
\end{matrix} \right].\end{equation}
\end{lemma}

\subsection{Some Preliminaries}
In this subsection, we give some preliminaries required in our following derivations.

\begin{lemma}[\cite{Matrix analysis}, Binet-Cauchy theorem] \label{Cauchy}
For an $n_1\times n_2$ matrix $M$ and an $n_2\times n_1$ matrix $N$, $n_1\le n_2$, it holds that
$$\det (MN) =\sum\nolimits_{J\subseteq [n_2], |J|=n_1} \det M_{J} \det (N^\intercal)_J.$$
\end{lemma}

\begin{lemma}[\cite{Matrix analysis}]\label{lemma3}
Let matrix $M=[M_1\ M_2]$ and $M_1^{\bot}$ be a basis of the left null space of $M_1$, i.e., a matrix consisting of a collection of row vectors that are linearly independent and span the left null space of $M_1$. Then, $M$ is of FRR, if and only if $M_1^{\bot}M_2$ is of FRR.
\end{lemma}

A {{{{matroid}} is a structure that abstracts the notion of linear independence in vector spaces \cite{Murota_Book}}}.
Given a finite set $E$  and a family ${{\cal I}}$ of subsets of $E$, the pair $(E, {\cal I})$ is a matroid if: (1)$\emptyset \in {{\cal I}}$; (2) if ${\cal I}_1 \in {\cal I}$ and ${\cal I}_2 \subseteq {{\cal I}_1}$, then ${\cal I}_2 \in {\cal I}$; (3) if ${\cal I}_1, {\cal I}_2 \in {\cal I}$ and $|{\cal I}_1|=|{\cal I}_2|+1$, then there is some $x\in {{\cal I}_1}\backslash {{\cal I}_2}$ satisfying ${\cal I}_2 \cup x \in {\cal I}$. {{In the above definition, $E$ is called the ground set, and a member of $\cal I$ is called an independent set. The rank of a matroid ${\cal M}=(E,{\cal I})$,  denoted by $\rho({\cal M})$, is defined as the maximum cardinality of its independent sets.}}
For a matrix $F$, a matroid can be defined  as ${\cal{M}}(F)=(E, {\cal I})$, where $E$ is the set of indices of columns of $F$, ${\cal I}$ is the collection of indices of columns of $F$ which are linearly independent, i.e., ${\cal I}=\{J\subseteq E: {{\rm rank}(F_{{J}})=|J|}\}$. This matroid is sometimes called a linear matroid. Given two matroids ${\cal{M}}_1$ and ${\cal M}_2$ over the same ground set, the matroid intersection ${\cal M}_1 \cap {\cal M}_2$ is defined as the collection of all common independent sets of ${\cal{M}}_1$ and ${\cal M}_2$. The cardinality of the largest independent set in ${\cal M}_1\cap {\cal M}_2$, which is also denoted by $\rho({\cal M}_1\cap {\cal M}_2)$ for notation simplicity, can be determined in polynomial time \cite{Murota_Book}. For two matroids ${\cal M}_1=({\cal I}_1, E)$ and ${\cal M}_2=({\cal I}_2, E)$, the matroid union ${\cal M}_1 \cup {\cal M}_2$ is a matroid $({\cal I}_3, E)$ such that any $X\in {\cal I}_3$ can be expressed as $X=Y\cup Z$ with $Y\in {\cal I}_1$ and $Z \in {\cal I}_2$. Determining the rank of the union of two linear matroids can also be done in polynomial time.

A set function $f: 2^{V}\rightarrow {\mathbb{R}}$ is {submodular} if for all sets $S_1 \subseteq S_2 \subseteq V$ and any element $s \in V\setminus {S_2}$, it holds that
$f(S_1\bigcup \{s\})-f(S_1) \ge f(S_2\bigcup \{s\})-f(S_2)$. This set function is nondecreasing if $f(S_2)\ge f(S_1)$ for all $S_1 \subseteq S_2 \subseteq V$. A well known fact is that the rank of a matrix is submodular, nondecreasing on any subset of its column vectors \cite{Murota_Book}.

 Given a directed graph (digraph) $\mathbb D$, let $V({\mathbb D})$ be its vertex set,  $E({\mathbb D})$ its edge set.
 A path from a vertex $v_1\in V({\mathbb D})$ to a vertex $v_n \in V({\mathbb D})$ is a sequence of edges $(v_1,v_2),(v_2,v_3),...,(v_{n-1},v_n)$ with no repeated end vertices, which is denoted by $\{v_1\rightarrow\cdots \rightarrow v_n\}$. If there is a path from $v_1$ to $v_2$, we say that $v_2$ is reachable from $v_1$. A path from a vertex to itself is called a cycle.  Length of a path is the number of its edges. {{A matching $M$ of a digraph ${\mathbb D}$ is a set of edges such that any of its element does not share a common start or end vertex with another element}}. The size of a matching $M$ is the number of edges contained in $M$. A digraph is strongly connected, if any two of its vertices are reachable from each other.  A strongly connected component (SCC) of a digraph is a subgraph that is strongly connected, and is maximal in the sense that no other edges or vertices can be included without breaking the property of being strongly connected. A clique is an undirected graph such that any two of its vertices are adjacent.

{{ Let $S\triangleq \{s_1,...,s_p\}$ be a set of algebraically independent variables, and $G(\lambda;S)$ be a TFM whose entries are rational functions of $\lambda$ with  coefficients being polynomials of $s_i|_{i=1}^{p}$ over $\mathbb{R}$. The generic rank of the TFM $G(\lambda;S)$ is the maximum rank that $G(\lambda;S)$ can achieve as a function of both $\lambda$ and $s_i|_{i=1}^p$. When $s_i|_{i=1}^{p}$ are fixed, $G(\lambda;S)$ is usually abbreviated as $G(\lambda)$, and the maximum rank that $G(\lambda)$ can achieve among all $\lambda$ is sometimes called the normal rank of $G(\lambda)$ \cite{Zhou_robust_book}.}} Zeros of a TFM have various definitions in the literature \cite{Kailath_1980}, \cite{Zhou_robust_book}. In this paper, given a $n_1\times n_2 $ dimensional TFM $G(\lambda)$ with a full row normal rank (FRNR), {{i.e., its normal rank equals its number of rows,}} we say that $\lambda_0 \in \mathbb{C}$ is a zero of $G(\lambda)$, if ${\rm rank}(G(\lambda_0))<n_1$. A TFM $G(\lambda;S)$ is said to have a zero depending on $S$, if it has a zero for arbitrarily fixed $s_i|_{i=1}^{p}$ in the corresponding parameter space, while the value of this zero is not independent of $s_i|_{i=1}^{p}$.

\section{Genericity of the Controllability of the NDS}
In this section, genericity is established for well-posedness and controllability of the NDS (\ref{SubsystemDynamics})-(\ref{SubsystemInteraction}).

Using the symbols $A_{xx}$, $A_{xv}$, $P$, etc. defined in Section II-B,  a lumped state-space representation can  be obtained for the NDS (\ref{SubsystemDynamics})-(\ref{SubsystemInteraction}), given as follows,
\begin{equation} \label{state space} \dot x(t) = Ax(t) + Bu(t) \end{equation}
where $x(t)={{\rm{{\bf{{col}}}}}}\{x_i(t)|_{i=1}^N\}$, $u(t)={{\rm{{\bf{{col}}}}}\{u_i(t)|_{i=1}^N\}}$, and
\begin{equation} \label{LFT parameterized} \begin{split}
[A\  B]&=[A_{xx}\  B_{xu}]{\rm{ + }}{A_{xv}} P {(I - {A_{zv}} P )^{ - 1}}[{{A_{zx}}}\  {{B_{zu}}}]\\
&={F_l}\left({\begin{pmat}[{.|.}]
{{A_{xx}}}&{{B_{xu}}}&{{A_{xv}}} \cr\-
{{A_{zx}}}&{{B_{zu}}}&{{A_{zv}}} \cr
\end{pmat}, P } \right),
\end{split}  \end{equation}
{{in which $F_l(\bullet,*)$ denotes the lower LFT operation defined in Section I.}}
As can be seen, the lumped matrix $[A\ B]$ also has a form of LFT parametrization. This has been already observed in \cite{zhou_2015}, and is a property of LFT expressions.

Let $s_1,...,s_k$ be the first principle parameters of $P$ in (\ref{scm}), and denote by $S=(s_1,...,s_k)$. Obviously, $k$ equals the number of nonzero entries in $P$. Let ${\bf{ P }}$ be the set of admissible matrices of $P$ parameterized by $S$ that make the NDS ${\bf{\Sigma}}$ well-posed, and ${\bf{S}} \subseteq {\mathbb{R}}^k$ the set of the corresponding admissible values of $S$. To clarify the dependence of matrices $A$ and $B$ on $S$ in (\ref{LFT parameterized}), in the following $A$ and $B$ are sometimes denoted by $A(S)$ and $B(S)$ respectively.

From the definition of well-posedness of a closed-loop system, it can be straightforwardly proven that the subsystem ${\bf{\Sigma}}_{i}$ is well-posed if and only if $\det (I{\rm{ - }}H^{(i)} P^{(i)}) \ne {\rm{0}}$, $i\in [N]$. Moreover, under the condition that each of its subsystems ${\bf{\Sigma}}_{i}|_{i=1}^N$, is well-posed, the NDS ${\bf{\Sigma}}$ is well-posed if and only if $\det (I{\rm{ - }}{A_{zv}} P) \ne {\rm{0}}$. Hence, ${\bf S}$ is open and dense in ${\mathbb{R}}^k$. With these results, it can be further shown that if there is a particular $ P$ such that the NDS ${\bf{\Sigma}}$ and each of its subsystems ${\bf{\Sigma}}_{i}|_{i=1}^N$, are well-posed, then for almost every $ P$ with the same structure, the corresponding system is also well-posed. {{ The essential ideas behind the proofs are that both $\det (I{\rm{ - }}{A_{zv}} P)$ and $\det(I{\rm{ - }}H^{(i)} P^{(i)})$ are polynomials of the nonzero elements of $P$, whose values vary continuously with the system first principle parameters. The details are omitted due to their obviousness.

A property of a parameter dependent system is called generic, if this property holds almost everywhere in the parameter space \cite{generic}. {{Obviously, $s_i|_{i=1}^k=0$ makes each subsystem ${\bf{\Sigma}}_i$ well-posed as well as the whole NDS. This implies that the NDS ${\bf{\Sigma}}$ is generically well-posed, irrespective of the fixed constants in (\ref{SubsystemDynamics}).}}}}

The following proposition reveals the genericity of controllability of the LFT parameterized NDS (\ref{SubsystemDynamics})-(\ref{SubsystemInteraction}).

\begin{proposition} \label{generic}
Controllability of the NDS  (\ref{SubsystemDynamics})-(\ref{SubsystemInteraction}) is a generic property.
\end{proposition}

{\emph{Proof:}} Consider the lumped representation of the NDS $\bf {\Sigma}$ (\ref{SubsystemDynamics})-(\ref{SubsystemInteraction}), which is given by (\ref{state space})-(\ref{LFT parameterized}). The controllability matrix of the system (\ref{state space}) is ${\cal C}(A(S),B(S))=[B(S)\  A(S)B(S)\ ...\  (A(S))^{M_x-1}B(S)]$. Suppose that there exists an $S_0\in {\bf S}$ such that the pair $(A(S_0),B(S_0))$ is controllable. This means that there is at least one set of $M_x$ linearly independent column vectors in ${\cal C}(A(S_0),B(S_0))$. Denote the matrix composed of these vectors by ${\cal C}(A(S_0),B(S_0))_J$, where $|J|=M_x$. It can be straightforwardly proven from (\ref{LFT parameterized}) that there exists two polynomials $d(S)$ and $n(S)$ of the first principle parameters $S$ such that $\det ({\cal C}(A(S),B(S))_J)=\frac{d(S)}{n(S)}$. Moreover, $n(S)$ can be expressed as $n(S)\triangleq \det(I- A_{zv} P)^{t}$ for some $t\in {\mathbb{N}}$. As ${\cal C}(A(S_0),B(S_0))_J$ is invertible, it further means that the polynomial $d(S)$ is not identically zero. Let $\bar {\bf S}=\{S\in {\bf S}: d(S)=0\}$. Then,  for any $S\in {\bf S}\backslash \bar {\bf S}$, the matrix ${\cal C}(A(S),B(S))_{J}$ is invertible and hence the pair $(A(S),B(S))$ is controllable, so is the NDS $\bf {\Sigma}$. Since the set $\bar {\bf S}$ has zero Lebesgue measure in ${\bf S}$, Proposition \ref{generic} follows. $\hfill\blacksquare$

{\normalsize{
Noting that both parametric inaccuracies and parametric variations are unavoidable in actual systems. The above results imply that rather than numerical values of the first principle parameters in each subsystem and the SCM $\Phi$, controllability of this NDS is mainly determined by how the system matrices are influenced by these parameters.}}

According to the well known PBH test, the matrix pair $(A,B)$ in (\ref{state space}) is controllable if and only if the matrix $[A-\lambda I\ B]$ is of FRR at each $\lambda\in \mathbb{C}$. In addition, a complex number $\lambda$ that makes $[A-\lambda I\ B]$ row rank deficient is called an uncontrollable mode \cite{Kailath_1980}.  As argued in \citep[Lemma 6.2]{Anderson_1982},
\citep[Lemma 2]{Rational_function}, \cite{Murota_Book}, and from Lemma \ref{lemma_4}, if System (\ref{state space})-(\ref{LFT parameterized}) parameterized by the first principle parameters $S$ is structurally uncontrollable (i.e., the associated system is uncontrollable for each fixed $S\in {\bf S}$), then the determinants of all the $(M_x+M_z)\times (M_x+M_z)$ dimensional submatrices of the MVP $M(\lambda)$ defined in (\ref{second_MVP}), expressed as polynomials of $s_i|_{i=1}^k$ and $\lambda$, share a common divisor. Moreover, this common divisor, if expressed as a polynomial of $\lambda$, has at least one root for each fixed $S\in {\bf S}$, which can be called the uncontrollable mode of System (\ref{state space})-(\ref{LFT parameterized}) with respect to (w.r.t.) parameters $S$.  Otherwise, this system is structurally controllable \cite{Murota_Book}.  To further characterize structural uncontrollability, uncontrollable modes are divided throughout this paper into the classes of fixed uncontrollable modes and parameter-dependent uncontrollable modes.

\begin{definition} \label{fixed uncontrollable mode} A fixed uncontrollable mode (FUM) of the NDS (\ref{SubsystemDynamics})-(\ref{SubsystemInteraction}) is a fixed $\lambda\in {\mathbb{C}}$ such that the matrix pair $[A(S)-\lambda I \ B(S)]$ is not of FRR for each fixed $S\in {\bf S}$. In other words, an FUM is an uncontrollable mode of the NDS (\ref{SubsystemDynamics})-(\ref{SubsystemInteraction}) that is independent of $S$.
\end{definition}

{{\begin{definition}  \label{varying uncontrollable mode} A parameter dependent uncontrollable mode (PDUM) of the NDS (\ref{SubsystemDynamics})-(\ref{SubsystemInteraction}) is an uncontrollable mode that is not independent of $S$; more precisely, it is an uncontrollable mode of the NDS (\ref{SubsystemDynamics})-(\ref{SubsystemInteraction}) for each fixed $S\in {\bf S}$, and across $S\in {\bf S}$ the set of its values does not belong to any subset of $\mathbb{C}$ with a finite cardinality. % any  set of a finite number of complex values.
\end{definition} }}

The above definitions are similar to the notion fixed mode suggested in \cite{Fixed Mode}.  According to the continuous dependence of roots of a polynomial on its coefficients, the PDUMs must piecewise continuously depend on $S$. Note that, whether or not a univariate polynomial with coefficients being polynomials of $s_i|_{i=1}^k$ over $\mathbb{R}$ has a root located in a fixed finite set of isolated complex values is a generic property in the parameter space of $S$ \cite{Murota_Book}. It can be shown that the existence of either an FUM or a PDUM is a generic property of the LFT parameterized system (\ref{state space})-(\ref{LFT parameterized}), so is of the NDS (\ref{SubsystemDynamics})-(\ref{SubsystemInteraction}) \citep[Section 6]{generic}.  In addition, let $\sigma_{f}$ denote the set of FUMs of System (\ref{state space})-(\ref{LFT parameterized}). From Definition \ref{fixed uncontrollable mode},
\begin{equation} \label{definition FUM}
\sigma_f=\bigcap\nolimits_{S\in {\bf S}}\{\lambda \in {\mathbb{C}}: [A(S)-\lambda I\ B(S)] {\rm \  is \  not \  of \  FRR}\}.
\end{equation}

Both FUM and PDUM play important roles in the following analysis about the structural controllability of the NDS ${\bf{\Sigma}}$.

\section{Structural controllability of LFT parameterized plants}
In order to derive computationally feasible conditions for the structural controllability of the NDS described by (\ref{SubsystemDynamics})-(\ref{SubsystemInteraction}), we at first establish a necessary and sufficient condition in this section for a general LFT parameterized matrix pair to be structurally controllable, that is, a matrix pair described by (\ref{LFT parameterized}) in which the matrices $A_{xx}$, $A_{xv}$, $A_{zx}$, $A_{zv}$, $B_{xu}$ and $B_{zu}$ are no longer block diagonal. Compared with \cite{Morse_1976}, \cite{Murota_Book} in which a linearly parameterized plant is adopted, the LFT parametrization adopted here is more natural in describing relations between system matrices of a plant and its first principle parameters and is able to represent a much wider class of plants. In addition, due to the existence of the inverse of the matrix $I- A_{zv} P$ in the parametrization, some mathematical challenges also arise.

For the sake of derivations, it is assumed that the matrices $A$ and $B$ in (\ref{LFT parameterized}) have the dimensions of $n\times n$ and $n \times q$ respectively. Moreover, the following assumption is adopted.\footnote{For an individual plant, in case that each first principle parameter has a rank-one coefficient matrix in $P$, Assumption \ref{ass2} can always be satisfied by introducing a transformation similar to (\ref{trans2LFT}) in Section V.} This assumption is eliminated in the next section for an NDS.
\begin{assumption}
$ P={ {\bf diag}} \{s_1,...,s_k \}$, where the parameters $s_1,...,s_k$ are algebraically independent.
\label{ass2}
\end{assumption}
\subsection{Existence of PDUMs via TFM Analysis}
This subsection is devoted to conditions of the existence of PDUMs via a structure analysis of some associated TFMs.

Define two TFMs respectively as
\[\begin{split}{G_{zv}}(\lambda ) &= {A_{zx}}{(\lambda I - {A_{xx}})^{ - 1}}{A_{xv}} + {A_{zv}},\\
{G_{zu}}\left(\lambda  \right) &= {A_{zx}}{(\lambda I - {A_{xx}})^{ - 1}}{B_{xu}} + {B_{zu}}.\end{split}\]
With the above definitions, the following Lemma \ref{lemma_proof_1} transforms the existence of PDUMs   into that of zeros of a TFM.

\begin{lemma} \label{lemma_proof_1}
The system (\ref{state space})-(\ref{LFT parameterized}) has a PDUM, if and only if the TFM
$\left[{{G_{zv}}(\lambda ) P  - I}\quad {{G_{zu}}(\lambda )}\right]$ has a zero depending on $S$.
\end{lemma}

\emph{Proof:} Note that when $\lambda \notin \sigma(A_{xx})$,
$[{{A_{zx}}{{(\lambda I - {A_{xx}})}^{ - 1}}}\ I]$ is a basis of the left null space of ${\bf{ col}}\{{\lambda I - {A_{xx}}}, { - {A_{zx}}}\}$.
Then by Lemma \ref{lemma3}, for a given $\lambda \notin \sigma(A_{xx})$, $M(\lambda)$  in (\ref{second_MVP}) is of FRR, if and only if $[
{{A_{zx}}{{(\lambda I - {A_{xx}})}^{ - 1}}}\quad I
 ]\left[ \begin{matrix}
{{B_{{xu}}}}&{ - {A_{xv}} P }\\
{{B_{zu}}}&{I - {A_{zv}} P }
\end{matrix} \right] = \left[
{{G_{zu}}(\lambda )}\ {I\! - \!{G_{zv}}(\lambda)}P \right]$ is of FRR.

Let ${\bf \Lambda}_c = \{\lambda \in \mathbb{C}: \lambda \notin \sigma(A_{xx})\}$.  Then a $\lambda$ belonging to ${\bf \Lambda}_c$ and varying with $S$ is a zero of $M(\lambda)$  of (\ref{second_MVP}), if and only if it is a zero of $[{{G_{zv}}(\lambda ) P  - I}\quad {{G_{zu}}(\lambda )}]$.
Notice that the set $\sigma(A_{xx})$ consists of only some isolated elements which do not affect the piecewise continuous dependence of the zeros of $[{{G_{zv}}(\lambda ) P  - I}\quad {{G_{zu}}(\lambda )}]$ on $S$. Hence, $M(\lambda)$ has a zero depending on $S$, if and only if $[{{G_{zv}}(\lambda) P  - I}\quad {{G_{zu}}(\lambda )}]$ has a zero depending on $S$. By Lemma \ref{lemma_4}, the result follows. $\hfill\blacksquare$

Recall that the TFMs $G_{zv}(\lambda)$ and $G_{zu}(\lambda)$ respectively have a dimension of $k\times k$ and a dimension of $k\times q$.
Construct a diagraph $\mathbb{L}=({\cal V}, {\cal E})$ associated with $[G_{zv}(\lambda)\  G_{zu}(\lambda)]$ as follows. Define the vertex set as ${\cal V}={\cal U}\bigcup {\cal Z}$ with ${\cal U}=\{u_1,...,u_q\}$ and ${\cal Z}=\{z_1,...,z_k\}$, while the edge set as ${\cal E}={\cal E}_{{\cal U}{\cal Z}}\bigcup {\cal E}_{{\cal Z}{\cal Z}}$ with ${{\cal E}_{{\cal U}{\cal Z}}} = \{ ({u_i},{z_j}):[{G_{zu}}{(\lambda )]_{ji}} \ne 0\} $, ${{\cal E}_{{\cal Z}{\cal Z}}} = \{ ({z_i},{z_j}):[{G_{zv}}{(\lambda )]_{ji}} \ne 0\}$.
The digraph constructed in this way is called the {{{auxiliary connection graph (ACG)}}} associated with $[G_{zv}(\lambda)\   G_{zu}(\lambda)]$. A vertex in ${\mathbb{L}}$ is input-reachable, if there exists a path from a vertex in $\cal U$ that ends at it. To clarify the dependence of an edge in ${\mathbb{L}}$ on the variable $\lambda$, edges of ${\cal E}_{{\cal Z}{\cal Z}}$ are classified by the following definitions.

\begin{definition} \label{edge_classify}
Given the ACG ${\mathbb L}$ associated with  $[G_{zv}(\lambda)\  G_{zu}(\lambda)]$, an edge $(z_i,z_j)\in {\cal E}_{{\cal Z}{\cal Z}}$ is a { $\lambda$-edge}, if $[G_{zv}(\lambda)]_{ji}$ does not take a value independent of the variable $\lambda$; otherwise $(z_i,z_j)\in {\cal E}_{{\cal Z}{\cal Z}}$ is a { constant-edge}. A cycle of $\mathbb{L}$ which contains at least one $\lambda$-edge is a {{$\lambda$-cycle}}. An { input-unreachable $\lambda$-edge} is a $\lambda$-edge such that neither its start vertex nor its end vertex is reachable from $\cal U$.  An input-unreachable $\lambda$-cycle (resp. input-unreachable SCC) is a $\lambda$-cycle (resp. SCC) with each of its vertices being input-unreachable.
\end{definition}

Based on the above definitions, the following proposition gives a graph theoretical necessary and sufficient condition for the existence of a PDUM of the system (\ref{state space})-(\ref{LFT parameterized}).

\begin{proposition} \label{proposition1}
Under Assumption \ref{ass2}, the system (\ref{state space})-(\ref{LFT parameterized}) has a PDUM, if and only if there exists an input-unreachable $\lambda$-cycle in the ACG $\mathbb{L}$.
\end{proposition}

{{ Proposition \ref{proposition1} implies that, under Assumption \ref{ass2},  the existence of a PDUM is only related to the structure of $[G_{zv}(\lambda)\  G_{zu}(\lambda)]$ and whether or not an entry of $G_{zv}(\lambda)$ depends on the variable $\lambda$, while  the interdependencies among different entries of $[G_{zv}(\lambda)\  G_{zu}(\lambda)]$ do not influence the obtained results.}}

To prove Proposition \ref{proposition1}, the following {Lemmas \ref{lemma_proof_mid_2}}-{\ref{lemma_proof_mid_4}} are required, whose proofs  are deferred to Appendix B. Lemma \ref{lemma_proof_mid_2} is an auxiliary result needed in the proofs of Lemmas \ref{lemma_proof_2} and \ref{lemma_proof_mid_3}. Lemma {\ref{lemma_proof_2}} presents a relation between the existence of zeros of a class of TFMs and that of a $\lambda$-cycle,  which is crucial in obtaining Proposition \ref{proposition1}. Lemmas {\ref{lemma_proof_mid_3}} and {\ref{lemma_proof_mid_4}} reveal a  relation between zeros of a class of TFMs and the input-reachability property of the associated ACG $\mathbb{L}$.

\begin{lemma}\label{lemma_proof_mid_2}
Given a digraph ${\mathbb G}$ with a vertex set $V=\{1,...,p\}$, $p\ge 2$, suppose that there are two distinct cycles $C_1$ and $C_2$ in ${\mathbb G}$ having an equal length $p$. Then, $\forall e \in { E}(C_1)$, there exists a cycle in ${\mathbb G}$  with length not exceeding $p-1$ that contains $e$.
\end{lemma}

\begin{lemma}\label{lemma_proof_2}
Under Assumption \ref{ass2}, let $\mathcal{J}=\{j_1,...,j_{n_s}\}\subseteq \{1,...,k\}$ and denote $(G_{zv}(\lambda))_{\mathcal{J},\mathcal{J}}$ by $G^\mathcal{J}_{zv}(\lambda)$, ${{\bf {\bf diag}}} \{s_{j_1},...,s_{j_{n_s}}\}$ by $P_{\mathcal{J}}$. Let $\mathbb{L}^{\mathcal{J}}$ be the subgraph of $\mathbb{L}$ induced by vertices $\{z_{j_1}$,...,$z_{j_{n_s}}\}$. Then, the TFM $G^\mathcal{J}_{zv}(\lambda)P_{\mathcal{J}}-I$ has a zero depending on parameters $s_{j_{i}}|_{i=1}^{n_s}$, if and only if $\mathbb{L}^{\mathcal{J}}$ has a $\lambda$-cycle.
\end{lemma}

\begin{remark}
 From the proof of Lemma \ref{lemma_proof_2}, it is clear that the necessity behind the introduction of $\lambda$-cycle is that, there are two types of nonzero entries in $G_{zv}(\lambda)$. As a result, the shortest cycle in the digraph ${{\mathbb D}}^{\cal J}$, which is defined in the proof of Lemma \ref{lemma_proof_2}, does not necessarily correspond to a nonvanishing term depending on $\lambda$ in ${\rm det} (G^\mathcal{J}_{zv}(\lambda)P_{\mathcal{J}}-I)$.
\end{remark}

\begin{lemma}\label{lemma_proof_mid_3} Under Assumption \ref{ass2}, if there does not exist an input-unreachable vertex $z\in {\cal Z}$ in the digraph $\mathbb{L}$ associated with $[G_{zv}(\lambda)\  G_{zu}(\lambda)]$, then the obtained matrix after deleting any column from  $[{{G_{zv}}(\lambda ) P  - I}\quad {{G_{zu}}(\lambda )}]$ is of full row generic rank (FRGR, {{i.e., its generic rank equals its number of rows}}).
\end{lemma}

\begin{lemma}\label{lemma_proof_mid_4}
Under Assumption \ref{ass2}, if the submatrix obtained by deleting any column from $[{{G_{zv}}(\lambda ) P  - I}\quad {{G_{zu}}(\lambda )}]$ is of FRGR,
then $[{{G_{zv}}(\lambda ) P  - I}\quad {{G_{zu}}(\lambda )}]$ does not have a zero depending on $S$.
\end{lemma}

{\emph{Proof of Proposition \ref{proposition1}:}} To prove the {{if}} part, suppose that there is an input-unreachable $\lambda$-cycle in the digraph $\mathbb{L}$. Let the vertex set of this input-unreachable $\lambda$-cycle be ${\cal Z}_s\triangleq\{z_{j_1},...,z_{j_{n_s}}\}\subseteq {\cal Z}$. Moreover, define a set $\cal J$ as ${\cal{J}}\triangleq \{j_1,...,j_{n_s}\}$. Since ${\cal Z}_s$ are unreachable from ${\cal U}$, it is clear that there exists a permutation matrix $Q$ such that \cite{generic} {\small{$Q{G_{zv}}(\lambda ){Q^{\intercal}} = \left[ \begin{matrix}
{{G^{\cal{J}}_{zv}}{{(\lambda )}}}&0\\
{{G^{[21]}_{zv}}{{(\lambda )}}}&{{G^{[22]}_{zv}}{{(\lambda )}}}
\end{matrix} \right], Q{G_{zu}}(\lambda ) = \left[ \begin{matrix}
0\\
{{G^{[2]}_{zu}}{{(\lambda )}}}
\end{matrix} \right],$}}  and therefore {\small{ \begin{equation} \label{lump}\begin{array}{l}
Q\left[\begin{matrix}
{{G_{zv}}(\lambda ) P  - I}&{{G_{zu}}(\lambda )}
\end{matrix}\right]\left[ \begin{matrix}
{{Q^{\intercal}}}&0\\
0&I
\end{matrix} \right] = \\
\begin{pmat}[{||}]
{{G^{\cal{J}}_{zv}}{{(\lambda )}}{P_{\cal{J}}} - I}&0&0\cr\-
{{G^{[21]}_{zv}}{{(\lambda )}}{P_{\cal{J}}}}&{{G^{[22]}_{zv}}{{(\lambda )}}{P_{[k]\backslash {\cal{J}}}}}-I&{{G^{[2]}_{zu}}{{(\lambda )}}}\cr
\end{pmat},
\end{array}\end{equation}}}where ${G^{[21]}_{zv}}{(\lambda )} \in F{(\lambda)^{(k-n_s) \times n_s}},{G^{[22]}_{zv}}{(\lambda )} \in F{(\lambda)^{(k-n_s) \times (k-n_s)}}$, ${G^{[2]}_{zu}}{(\lambda )} \in F{(\lambda)^{(k - {n_s}) \times q}}$, $G^{\cal J}_{zv}(\lambda)$ and $P_{\cal J}$ are defined in Lemma \ref{lemma_proof_2}. By Lemma \ref{lemma_proof_2}, $G^{\cal{J}}_{zv}{{(\lambda )}}{P_{\cal J}} - I$ has a zero depending on $s_{j_{i}}|_{i=1}^{n_s}$, so does $[{{G_{zv}}(\lambda ) P  - I}\quad {{G_{zu}}(\lambda )}]$  noting that according to (\ref{lump}) the associated ranks remain invariant. By Lemma \ref{lemma_proof_1}, this means that the system (\ref{state space})-(\ref{LFT parameterized}) has a PDUM.

To prove the {only if} part, first assume that every vertex in ${\cal Z}$ is input-reachable in $\mathbb{L}$. Then from Lemmas \ref{lemma_proof_mid_3}-\ref{lemma_proof_mid_4}, we have that $[{{G_{zv}}(\lambda ) P  - I}\quad{{G_{zu}}(\lambda )}]$ does not have a zero depending on $S$.  Now suppose that there is an input-unreachable vertex in $\mathbb{L}$, and denote the set of all input-unreachable vertices by ${ \cal \bar Z}_s \subseteq {\cal Z}$. Suppose that there does not exist a $\lambda$-cycle in the subgraph of $\mathbb{L}$ induced by ${\cal \bar Z}_s$ (denoted by ${\cal L}^{{\cal \bar Z}_s}$).  Following (\ref{lump}), there is a permutation matrix $\bar Q$ such that {\small{
\begin{equation} \label{matrix_transform_2}  \begin{array}{l}
{\bar Q}\left[\begin{matrix}
{{G_{zv}}(\lambda ) P  - I}&{{G_{zu}}(\lambda )}
\end{matrix}\right]\left[ \begin{matrix}
{\bar{Q}^{\intercal}}&0\\
0&I
\end{matrix} \right] = \\
\begin{pmat}[{||}]
{\bar G_{zv}^{[11]}{{(\lambda )}}\bar P_1 - I}&0&0\cr\-
{\bar G_{zv}^{[21]}{{(\lambda )}}{\bar P_1}}&{\bar G_{zv}^{[22]}{{(\lambda )}}\bar P_2-I}&{\bar G_{zu}^{[2]}{{(\lambda )}}}\cr
\end{pmat}
\end{array},\end{equation}}}where $\bar G_{zv}^{[11]}{(\lambda )} \in F{(\lambda)^{|{{\cal \bar Z}_s}| \times |{{\cal \bar Z}_s}|}}$, $\bar G_{zv}^{[21]}{(\lambda )} \in F{(\lambda)^{|{\cal Z}\backslash {{\cal \bar Z}_s}| \times |{{\cal \bar Z}_s}|}}$, $\bar G_{zv}^{[22]}{(\lambda )} \in F{(\lambda)^{|{\cal Z}\backslash {{\cal \bar Z}_s}| \times |{\cal Z}\backslash {{\cal \bar Z}_s}|}}$, $\bar G_{zu}^{[2]}{(\lambda )} \in F{(\lambda)^{|{\cal Z}\backslash {{\cal \bar Z}_s}| \times q}}$ and ${{\bf {\bf diag}}}\{\bar P_1,\bar P_2\}\triangleq \bar Q  P {\bar {{Q}}^{\intercal}}$. Since there is not a $\lambda$-cycle in the digraph ${\cal L}^{{\cal \bar Z}_s}$ associated with ${\bar G^{[11]}_{zv}{{(\lambda )}}}$, it can be obtained  from Lemma \ref{lemma_proof_2} that
${\bar G^{[11]}_{zv}{{(\lambda )}}\bar P_1 - I}$ does not have a zero depending on $S$. Meanwhile, noting that ${\cal Z} \backslash {\cal \bar Z}_s$ is the set of all input-reachable vertices,  Lemmas \ref{lemma_proof_mid_3}-\ref{lemma_proof_mid_4} imply that $[\bar G_{zv}^{[22]}{{(\lambda )}}\bar P_2-I\quad {\bar G_{zu}^{[2]}{{(\lambda )}}}]$ does not have a zero depending on $S$.  Notice that the right-hand side matrix of (\ref{matrix_transform_2}) has a block triangular structure. Hence, it does not have a zero depending on $S$ either. This indicates the nonexistence of a PDUM in the original system (\ref{state space})-(\ref{LFT parameterized}) by Lemma \ref{lemma_proof_1}.  $\hfill\blacksquare$
\subsection{Existence of FUMs}
From the definition of an FUM, since $0\in {\bf S}$, it suffices to see that the set of FUMs $\sigma_f$ of System (\ref{state space})-(\ref{LFT parameterized}) satisfies
 $\sigma_f\subseteq \sigma(A_{xx})$.{\footnote{This inclusion relation is still valid even when the first principle parameters $s_i|_{i=1}^k$ can only take values from some continuous intervals not including $0$. This is because an FUM must be a zero of certain polynomial of $\lambda$ involving $s_i|_{i=1}^k$, which can also be regarded as a polynomial of each $s_i$ by substituting the value of this FUM and fixing the rest of $s_i|_{i=1}^k$. As such polynomial equals zero for infinitely many values of $s_i$ from Definition \ref{fixed uncontrollable mode}, it must be a zero polynomial of $s_i$. Hence, fixing $s_i|_{i=1}^k$ to be an arbitrary real value satisfying the well-posedness condition, the set of uncontrollable modes of the obtained system must contain $\sigma_f$.}} Hence, the existence of an FUM is equivalent to that of a $\lambda_0 \in \sigma(A_{xx})$, such that $M(\lambda_0)$  in (\ref{second_MVP}) is not of FRR for {{every}} ${ P} \in {\bf{ P }}$.

The following proposition presents an equivalent condition for the above statements in terms of matroid union.

\begin{proposition} \label{Proposition 2}
Under Assumption \ref{ass2}, the system (\ref{state space})-(\ref{LFT parameterized}) has an FUM, if and only if there exists a $\lambda_0 \in \sigma(A_{xx})$, so that $\rho({{\cal M}_1} \cup {{\cal M}_2(\lambda_0)}) < 2k + n$, in which matroids
${{\cal M}_1} = {\cal {\cal M}}([{0_{k\times n}}\  {I_k}\  I_k])$, ${{\cal M}_2(\lambda_0)} = {\cal {\cal M}}\left({\left[ \begin{matrix}
{\lambda_0 I - {A_{xx}}}&{{B_{xu}}}&{{A_{xv}}}\\
{{-A_{zx}}}&{{B_{zu}}}&{{A_{zv}}}\\
0&0&I
\end{matrix} \right]^{\intercal}}\right)$.
\end{proposition}

{\emph{Proof:}}  Define an MVP $H(\lambda)$ as
\begin{equation}\label{basic_MVP} H(\lambda)=\left[ \begin{matrix}
{\lambda I - {A_{xx}}}&{{B_{xu}}}&{{A_{xv}}}&0\\
{ - {A_{zx}}}&{{B_{zu}}}&{{A_{zv}}}&I\\
0&0&I&{{ P}}
\end{matrix} \right],\end{equation} where $P\in {\bf P}$.   Using Schur complement \cite{Matrix analysis}, it can be validated that for any fixed $\lambda \in \mathbb{C}$ and $P$, $M(\lambda)$ in (\ref{second_MVP}) is of FRR if and only if $H(\lambda)$ is of FRR.
Now substitute $P={\rm{{\bf{{\bf diag}}}}}\{s_i|_{i=1}^k\}$ into $H(\lambda_0)$. Note that
\begin{equation} \label{transform}
\begin{array}{l}
H({\lambda _0})\left[ \begin{matrix}
I_{n+q+k}&0\\
0&{{{\bf {\bf diag}}}\{ {t_i}\left| {_{i = 1}^k} \right.\} }
\end{matrix} \right] = \\
\left[ \begin{matrix}
{{\lambda _0}I - {A_{xx}}}&{{B_{xu}}}&{{A_{xv}}}&0\\
{{-A_{zx}}}&{{B_{zu}}}&{{A_{zv}}}&{{{\bf {\bf diag}}}\{ {t_i}\left| {_{i = 1}^k} \right.\} }\\
0&0&I&{{{\bf {\bf diag}}}\{ {t_i}{s_i}\left| {_{i = 1}^k} \right.\} }
\end{matrix} \right]
\end{array}\end{equation} where auxiliary variables $t_1,...,t_k$ are nonzero, and $\{t_1,...,t_k,s_1,...,s_k\}$ are algebraically independent. It suffices to see that the rank of $H(\lambda_0)$ equals that of the right-hand side of (\ref{transform}).
 Notice that $\{s_1,...,s_k, s_1t_1,...,s_kt_k\}$ are also algebraically independent. Based on Theorem 4.2.3 of \cite{Murota_Book}, for a matrix $F={\bf{{col}}}\{Q,T\}$ where $Q$ is a numerical matrix over $\mathbb{C}$ and $T$ is a structured matrix with its nonzero entries being algebraically independent, it holds that $${\cal M}(F)={\cal M}(Q)\cup {\cal M}(T).$$
 Using the above result on the transpose of the right-hand side matrix of (\ref{transform}), the required result follows from the fact that ${\cal{M}}([0, {{\bf {\bf diag}}}\{t_i\left| {_{i = 1}^k} \right.\}, {{\bf {\bf diag}}}\{{t_i}{s_i}\left| {_{i = 1}^k} \right.]\})={\cal{M}}([0\ I_k\ I_k])$. $\hfill\blacksquare$

~\\
%{{{\subsection{Necessary and Sufficient Condition for Structural Controllability}}}}

Combining Propositions {\ref{proposition1}} and {\ref{Proposition 2}}, we have the following theorem.

\begin{theorem}\label{theorem1}
Under Assumption \ref{ass2}, the LFT parameterized plant (\ref{state space})-(\ref{LFT parameterized}) is structurally controllable, if and only if the following two conditions hold simultaneously.

(i) There is no input-unreachable $\lambda$-cycle in the ACG $\mathbb{L}$;

(ii) For each $\lambda_0 \in \sigma(A_{xx})$, $\rho({\cal M}_1 \cup {{\cal M}_2(\lambda_0)}) = 2k + n$.
\end{theorem}

{\emph{Proof:}} From Propositions {\ref{proposition1}} and {\ref{Proposition 2}}, the necessity is obvious. For sufficiency, the simultaneous satisfaction of Conditions (i) and (ii) implies that the plant (\ref{state space})-(\ref{LFT parameterized}) cannot be structurally uncontrollable. In other words, it is structurally controllable.  $\hfill\blacksquare$

Checking whether an ACG $\mathbb{L}$ associated with a given system has an input-unreachable $\lambda$-cycle can be done efficiently using  SCC decompositions with a complexity of $O(|V({\mathbb L})|+|E({\mathbb L})|)$ \cite{DB_West_graph}, which is briefly as follows:
%\begin{itemize}

(1) Do SCC decompositions on ${\mathbb L}$, and find all the input-unreachable SCCs;

(2) If there is a $\lambda$-edge in at least one input-unreachable SCC, there will be an input-unreachable $\lambda$-cycle in $\mathbb{L}$. Otherwise there will be none.
%\end{itemize}

The rationale of the above procedure lies in that, all vertices of a $\lambda$-cycle must belong to the same SCC, as they are reachable from each other. Besides, if there is a $\lambda$-edge in an SCC, then a $\lambda$-cycle exists, noting that the start and end vertices of this $\lambda$-edge are mutually reachable.  The matroid union condition (ii) of Theorem {\ref{theorem1}} enables efficient verification with a polynomial time complexity.

In addition, Theorem {\ref{theorem1}} provides some other information for a structurally uncontrollable system. Particularly, if a plant has a PDUM, there must exist some uncontrollable modes that depend on the first principle parameters associated with the shortest input-unreachable $\lambda$-cycle in the ACG $\mathbb{L}$, which is indicated in the proof of Lemma \ref{lemma_proof_2}. Moreover, any $\lambda_0$ that fails to satisfy Condition (ii) of Theorem \ref{theorem1} is an FUM of the plant.

{\begin{remark} \label{remark_difference} A necessary and sufficient condition of structural controllability for a linear parameterized plant was first given in \cite{Morse_1976}. The condition there {{has exponential time complexity}}, and was derived from the decentralized stabilization theory \cite{Morse_decentralized}, which cannot be directly extended to LFT parameterized plants due to the existence of the matrix inverse $(I-A_{zv} P)^{-1}$. The descriptor system approach of \cite{Murota_SIAM} is under the nondimensionality assumption defined therein, which introduces some restrictions on the fixed nonzero constants of the associated matrices. Our derivations for the LFT parameterized plant are mainly based on structure analysis of the associated TFMs, and do not impose any assumptions on the associated constant matrices. A concept ``$\lambda$-cycle'' is introduced, which is important in establishing Theorem \ref{theorem1}. {{Moreover, these results have some physical interpretations when applied to an NDS as shown in the next section.}} It is worthwhile to mention that the approach suggested here is also capable of handling some situations beyond the rank-one case. Due to space considerations, these extensions are not discussed here.
\end{remark}

\section{A Network based necessary and sufficient condition for structural controllability}

In this section, results in Section IV are extended to the NDS (\ref{SubsystemDynamics})-(\ref{SubsystemInteraction}) by {{taking the system matrix structures into account}}. This leads not only to a more computationally efficient method for structural controllability verification of an NDS, but also to some deeper insights on how subsystem dynamics and subsystem interconnection topology influence controllability of the whole NDS. {{Please note that Assumption \ref{ass2} is not needed in this section.}}

For Subsystem ${\bf{\Sigma}}_i$ in (\ref{SubsystemDynamics}) with its augmented system matrices given in (\ref{lft_sub}), define TFMs $G_{zv}^{(i)}(\lambda)$ and $G_{zu}^{(i)}(\lambda)$ respectively as
\[ \begin{split} G_{zv}^{(i)}(\lambda ) =  A_{zx}^{(i)}{(\lambda I -  A_{xx}^{(i)})^{ - 1}} A_{xv}^{(i)} +  A_{zv}^{(i)},\\
G_{zu}^{(i)}(\lambda ) =  A_{zx}^{(i)}{(\lambda I -  A_{xx}^{(i)})^{ - 1}} B_{xu}^{(i)} +  B_{zu}^{(i)}.\end{split}\]Construct the ACG $\mathbb{T}_i$ associated with $[G_{zv}^{(i)}(\lambda)\  G_{zu}^{(i)}(\lambda )]$ as
$\mathbb{T}_i=({\cal U}_i\cup {\cal V}_i \cup {\cal Z}_i, {\cal E}_{{\cal U}_i{\cal Z}_i}\cup {\cal E}_{{\cal V}_i{\cal Z}_i})$, where ${\cal U}_i=\{u_{i1},...,u_{im_{ui}}\}$, ${\cal V}_i=\{v_{i1},...,v_{im_{vi}}\}$ and ${\cal Z}_i=\{z_{i1},...,z_{im_{zi}}\}$ represent the vertex set of external inputs, internal inputs and internal outputs of $\tilde {\bf \Sigma}_i$ respectively; ${\cal E}_{{\cal V}_i{\cal Z}_i}=\{(v_{ip},z_{iq}): [G^{(i)}_{zv}(\lambda)]_{qp}\ne 0, p\in[m_{vi}], q\in[m_{zi}]\}$, ${\cal E}_{{\cal U}_i{\cal Z}_i}=\{(u_{ip},z_{iq}): [G^{(i)}_{zu}(\lambda)]_{qp}\ne 0, p\in[m_{ui}], q\in[m_{zi}]\}$.   Construct a digraph $\mathbb{T}_{{\bf \Sigma}}=({\cal V}_{{\bf \Sigma}},{\cal E}_{{\bf \Sigma}})$ by connecting these $N$ ACGs $\mathbb{T}_i$ from ${\cal Z}_i$ to ${\cal V}_j$  through the edge set ${\cal E}_{P}=\{(z_{ip},v_{jq}): [P_{j,i}]_{qp} \ne 0, i,j \in[N], p\in[m_{zi}], q\in[m_{vj}]\}$ where $  P_{i,j}$ is the $(i,j)$-th block of $  P$ in (\ref{scm}). That is, ${\cal V}_{{\bf \Sigma}}=\bigcup\nolimits_{i =1}^{N} {\left( {{{\cal U}_i}\cup {{{\cal Z}_i}\cup {{{\cal V}_i}} } } \right)}$, ${\cal E}_{{\bf \Sigma}}={{\cal E}_{  P} }\bigcup \{\bigcup\nolimits_{i =1}^N {\left({{{\cal E}_{{{\cal U}_i}{{\cal Z}_i}}}\cup {{{\cal E}_{{{\cal V}_i}{{\cal Z}_i}}}}} \right)}\}$. Afterwards, this digraph is called the networked-ACG, and abbreviated as n-ACG. {{Recall that an entry of $G_{zv}^{(i)}(\lambda)$ which depends on the variable $\lambda$ corresponds to a $\lambda$-edge in ${\mathbb{T}}_i$, and definitions of other graph elements in $\mathbb{T}_{{\bf \Sigma}}$, like the $\lambda$-cycle, inherit Definition \ref{edge_classify}.}}   An illustration of n-ACGs can be found in Fig. \ref{n-ACG} of Section VII. { Notice that $G^{(i)}_{zv}(\lambda)$ and $G^{(i)}_{zu}(\lambda)$ are respectively the TFMs from the internal inputs $\tilde v_i(t)$ and the external inputs $u_i(t)$ to the internal outputs $\tilde z_i(t)$ of the augmented system of Subsystem ${\bf \Sigma}_i$, and $P$ has the same structure as the SCM $\Phi$ except for its diagonal blocks. It seems safe to declare that $\mathbb{T}_{{\bf \Sigma}}$ intuitively reflects the {information flows over the NDS (\ref{SubsystemDynamics})-(\ref{SubsystemInteraction})}. More specifically, the existence of a path between two vertices in $\mathbb{T}_{{\bf \Sigma}}$ means that, the associated internal input/output variable of one subsystem can receive signal from the other associated one in the NDS.}

For the NDS (\ref{SubsystemDynamics})-(\ref{SubsystemInteraction}) and $P$ in (\ref{scm}), to meet Assumption \ref{ass2}, consider the following transformation on its lumped representation (\ref{LFT parameterized})
\begin{equation} \label{trans2LFT}{\small{
 \begin{split}
{F_l}\left(\!\begin{pmat}[{.|.}]
{{{ A}_{xx}}}&{{{ B}_{xu}}}&{{{ A}_{xv}}}\cr\-
{{{ A}_{zx}}}&{{{ B}_{zu}}}&{{{ A}_{zv}}}\cr
\end{pmat},  P  \!\!\right)
\!\!=\!\! {F_l}\left(\!  \begin{pmat}[{.|.}]
{{{ A}_{xx}}}&{{{ B}_{xu}}}&{{{ A}_{xv}}U}\cr \-
{V{{ A}_{zx}}}&{V{{ B}_{zu}}}&{V{{ A}_{zv}}U}\cr
\end{pmat},P_d \!\!\right)
 \end{split}}}
\end{equation}
where $P_d \triangleq {\rm{{\bf{{\bf {\bf diag}}}}}}\{s_i|_{i=1}^{k}\}$, $U$ and $V$ are some constant matrices satisfying $ P= UP_dV$. Here, each column of the matrix $U$ and each row of the matrix $V$ is a {{standard}} unit basis vector of an Euclidean space. Under the condition that all the first principle parameters of the NDS are algebraically independent, it can be simply shown that this decomposition always exists.

 For the lumped representation in the right-hand side of (\ref{trans2LFT}), define TFMs $G^{{\bf \Sigma}}_{zv}(\lambda)$ and $G^{{\bf \Sigma}}_{zu}(\lambda)$ respectively as $G_{zu}^{\bf \Sigma} \left( \lambda  \right)\!\! \triangleq \!\! V{A_{zx}}{(\lambda I - {A_{xx}})^{ - 1}}{B_{xu}} + V{B_{zu}}$, $G_{zv}^{\bf \Sigma} (\lambda )\! \triangleq\!\! V{A_{zx}}{(\lambda I - {A_{xx}})^{ - 1}}{A_{xv}}U + V{A_{zv}}U\!\!$. Then obviously $G^{{\bf \Sigma}}_{zv}(\lambda)\!=\! VG_{zu}\left( \lambda  \right)$, $G_{zv}^{\bf \Sigma}(\lambda )\!= \! VG_{zv}(\lambda )U$,
where $G_{zu}(\lambda ){\rm{ = {\bf{{\bf {\bf diag}}}}}}\{ G_{zu}^{(i)}(\lambda )|_{i = 1}^N\}$ and $G_{zv}{\rm{ = {\bf{{\bf {\bf diag}}}}}}\{ G_{zv}^{(i)}(\lambda )|_{i = 1}^N\}$. Denote the ACG associated with $[G^{{\bf \Sigma}}_{zv}(\lambda)\  G^{{\bf \Sigma}}_{zu}(\lambda)]$ by $\mathbb{L}_{{\bf \Sigma}}$. {{We have the following result, which establishes a relation between $\mathbb{L}_{\bf \Sigma}$ and $\mathbb{T}_{\bf \Sigma}$ w.r.t. the existence of an input-unreachable $\lambda$-cycle, thus removing Assumption \ref{ass2}. Its proof is deferred to Appendix B, where the essential idea is to build a relation between graph connectivity and multiplications of the associated matrices. A similar issue is addressed in \cite{Tsopelakos bilinear} for bilinear systems by using unit matrices.}}

\begin{proposition}\label{proposition_3}
There exists an input-unreachable $\lambda$-cycle in $\mathbb{L}_{{\bf \Sigma}}$, if and only if there exists an input-unreachable $\lambda$-cycle in the n-ACG $\mathbb{T}_{{\bf \Sigma}}$. In other words, the NDS (\ref{SubsystemDynamics})-(\ref{SubsystemInteraction}) does not have a PDUM, if and only if there is no input-unreachable $\lambda$-cycle in $\mathbb{T}_{{\bf \Sigma}}$.
\end{proposition}

To derive a computationally efficient criterion for verification of an FUM for NDSs, we adopt some ideas  similar to those in \cite{zhou_2015} and \cite{Y_Zhang_2016}. For each subsystem,  define an MVP ${M_i}(\lambda ) \triangleq \left[\begin{matrix}
{\lambda I - A_{xx}^{(i)}}&{B_{xu}^{(i)}}\\
{ - A_{zx}^{(i)}}&{B_{zu}^{(i)}}
\end{matrix} \right]$.
Suppose that there are $m$ distinct values in $\bigcup\nolimits_{i\in[N]} \sigma(A^{(i)}_{xx})$, and denote the set constituted from them by ${\bf \Lambda}= \{\lambda_1,...,\lambda_m\}$.
For $i\in [m]$,  $j\in[N]$,  let $[{T^{(j)}_{i}}\  {Z^{(j)}_{i}}]$ be a matrix constituted by a basis of the left null space of $M_{j}(\lambda_i)$, and assume that $T^{(j)}_{i}$ and $Z^{(j)}_{i}$ have dimensions ${m_{rij}\times m_{xj}}$ and ${m_{rij}\times m_{zj}}$ respectively.  Notice that $m_{rij}=0$ if $\lambda_i$ is not a zero of this $M_j(\lambda)$ which is of FRNR.
Moreover, let $M_{ri}=\sum\nolimits_{j=1}^{N}m_{rij}$, and ${Y^{(j)}_{i}} = {T^{(j)}_{i}}A_{xv}^{(j)} + {Z^{(j)}_{i}}A_{zv}^{(j)}$. Construct matrices\footnote{In case that $m_{rij}=0$, the number of rows of $T_i$ does not increase. The only things needing to do are to put $m_{xj}$ zeros to the corresponding column entries. Similar remarks are applicable to $Y_i$ and $Z_i$.}
$T_i={{\rm{{\bf{{\bf {\bf diag}}}}}}}\{ T^{(j)}_{i}|_{j=1}^N\}$, $Y_i= {\rm{{\bf{{\bf {\bf diag}}}}}}\{ Y^{(j)}_{i}|_{j=1}^N\}$, $Z_i = {\rm{{\bf{{\bf {\bf diag}}}}}}\{ {Z^{(j)}_{i}|_{j=1}^N}\}$. It is clear that $[T_i\  Z_i]$ is a basis of the left null space of ${\bf{{col}}}\{[\lambda_i I-A_{xx}\  B_{xu}],[-A_{zx}\  B_{zu}]\}$. By Lemma \ref{lemma3}, $M(\lambda_i)$ defined in (\ref{second_MVP}) is of {{FRGR}}, if and only if
$$\left[{{T_i}}\ {{Z_i}}\right]\left[ \begin{matrix}
{ - {A_{xv}} P }\\
{I - {A_{xv}} P }
\end{matrix} \right]$$ is of {{FRGR}},
which is further equivalent to that
$[ {{Y_i}}\ Z_i ]{\bf col}\{ P , I_{M_z}\}$ is of {{FRGR}}. Hence, we have the following proposition.

\begin{proposition} \label{proposition 4} The NDS {$\bf{\Sigma}$} (\ref{SubsystemDynamics})-(\ref{SubsystemInteraction}) does not have an FUM, if and only if for each $i\in[m]$, ${\rho}({\cal M}({Q_1}) \cap {\cal M}({Q_{2i}})) = {M_{ri}}$, where
${Q_1}{\rm{=}}[{{{ P ^{\intercal}}}}\ {{I_{{M_z}}}}],{Q_{2i}}{\rm{=}}[{Y_i}\ {Z_i}]$.
\end{proposition}

{\emph{Proof:}}
Notice that for a set of algebraically independent variables $\{s_1,...,s_k,t_1,...,t_{M_z}\}$ where each $t_i$ is nonzero,
\begin{equation}\label{Zhuanhuan}{\small{
\left[{{Y_i}}\ {{Z_i}}\right]\left[ \begin{matrix}
 P \\
I
\end{matrix} \right]{\rm{{\bf{{\bf {\bf diag}}}}}}\{ {t_1},...,{t_{{M_z}}}\}=\left[
{{Y_i}}\ {{Z_i}}\right]\left[\begin{matrix}
{ P {\rm{{\bf{{\bf {\bf diag}}}}}}\{ {t_1},...,{t_{{M_z}}}\} }\\
{{\rm{{\bf{{\bf {\bf diag}}}}}}\{ {t_1},...,{t_{{M_z}}}\} }
\end{matrix} \right].}}
\end{equation}It means that the generic rank of $[Y_i\  Z_i]{{\bf{{col}}}\{ P,I_{M_z}\}}$  equals that of the right-hand side of (\ref{Zhuanhuan}). Denote ${\rm{{\bf{{\bf {\bf diag}}}}}}\{ {t_1},...,{t_{{M_z}}}\}$ by $\Pi$. Note that the nonzero entries of ${\bf col}\{ P \Pi, \Pi\}$ are algebraically independent. Hence, ${\cal{M}}([\Pi P^{\intercal}\ \Pi])={\cal M}([P^{\intercal}\ I_{M_z}])$.  From Lemma \ref{Cauchy}, it can be directly shown that the right-hand side of (\ref{Zhuanhuan}) is of FRGR, if and only if there exists a $J\subseteq [M_v+M_z]$ with $|J|=M_{ri}$ such that $Q_{2i,J}$ is of full column rank, while $[\Pi P^{\intercal}\ \Pi]_J$ is of FRGR. By the definition of matroid intersection, this is further equivalent to that $\rho({\cal M}({Q_1}) \cap {\cal M}({Q_{2i}})) = {M_{ri}}$.
 $\hfill\blacksquare$

As a result of Propositions \ref{proposition_3} and \ref{proposition 4}, the structural controllability verification of the NDS (\ref{SubsystemDynamics})-(\ref{SubsystemInteraction}) can be implemented efficiently. In Proposition \ref{proposition_3}, the arithmetic operations are imposed on each subsystem in constructing the n-ACG $\mathbb{T}_{{\bf \Sigma}}$.
The existence of an input-unreachable $\lambda$-cycle can be checked using SCC decompositions in time complexity $O(M^2_x)$. In Proposition \ref{proposition 4}, parameters constituting $Y_i$ and $Z_i$ can be constructed within each subsystem because of their block diagonal structures.  When the matroid intersection algorithm is utilized, the only operation that involves arithmetic calculations is the independence oracle call, i.e., verifying whether columns of $Q_{1{{J}}}$ or $Q_{{2i,{J}}}$ for ${{J}}\subseteq [M_v+M_z]$ are linearly independent \cite{Murota_Book}. Denote the time of such an operation by $\tau$. Then, the conditions of Proposition \ref{proposition 4} can be tested in polynomial time with complexity $O(\tau M_x(M_z+M_v)^{2.5})$ \cite{Murota_Book}.  Note that $Q_{2i}$ composes of two column blocks with each one being block diagonal, and $Q_{1}$ can be seen structured.  Hence, verifying the conditions of both Propositions \ref{proposition_3} and \ref{proposition 4}  involves arithmetic operations  only within each subsystem, and graphical operations on the network topology. This property makes our approach attractive when dealing with large-scale systems.

Propositions \ref{proposition_3} and \ref{proposition 4} also illustrate how subsystem input-output relations and subsystem uncontrollable modes, together with the subsystem interconnection topology, jointly influence network controllability. {{Note that the n-ACG ${\mathbb{T}}_{\bf{\Sigma}}$ somehow reflects the information flows over the NDS, and a $\lambda$-edge means that this channel contains some state information of one associated subsystem. Proposition \ref{proposition_3} means that, to guarantee the absence of a PDUM, it is necessary and sufficient to make sure that there does not exist any state involved closed loop which cannot receive signal from the external inputs.}} Moreover, it can be validated that, if $\lambda_i$ is an uncontrollable mode shared by some subsystems $\tilde {\bf \Sigma}_{j_s}|_{s=1}^{p}$ (i.e., $[A^{(j_s)}_{xx}-\lambda_i I\ B^{(j_s)}_{xu}]$ is of row rank deficient), then $Z^{(j_s)}_{i}$ is not of FRR, $s\in [p]$. From the property of matroid intersection, Proposition \ref{proposition 4} implies that, to guarantee that $\lambda_i$ is not an FUM of the NDS,  at least $\sum\nolimits_{s=1}^p(m_{rij_s}\!\!-\!{\rm rank} (Z^{(j_s)}_{i}))$, which is equal to $\sum\nolimits_{s=1}^p(m_{xj_s}\!\!-\!{\rm rank}([\lambda_i I-A^{(j_s)}_{xx}\ B^{(j_s)}_{xu}]))$(\cite{Matrix analysis}),  vertex-disjoint subsystem interconnection links (corresponding to the nonzero entries in $P$), including possible self-loops, should be injected to these augmented subsystems $\tilde {\bf \Sigma}_{j_s}|_{s=1}^p$. These observations are helpful in constructing a controllable NDS.

Finally, Propositions \ref{proposition_3} and \ref{proposition 4} can be combined
to derive a simpler criterion for structural controllability of the NDS.

\begin{corollary} \label{corollary_add} The NDS $\bf{\Sigma}$ (\ref{SubsystemDynamics})-(\ref{SubsystemInteraction}) is structurally controllable, if and only if

(1) There is no input-unreachable {{$\lambda$-edge}} in the n-AGC ${\mathbb{T}}_{{\bf \Sigma}}$;

(2) For each $i\in[m]$, ${\rho}({\cal M}({Q_1}) \cap {\cal M}({Q_{2i}})) = {M_{ri}}$, where
${Q_1}$ and $Q_{2i}$ are defined in Proposition \ref{proposition 4}.

\end{corollary}

{\emph{Proof:}} The satisfaction of Condition (1) of Corollary {\ref{corollary_add}} certainly leads to that of the condition of Proposition \ref{proposition_3}. Hence, the {if} part is obvious. To show the {only if} part, suppose that there exists an input-unreachable $\lambda$-edge $(v_{ip}, z_{iq})$ which is not contained by any cycle in the n-ACG ${\mathbb{T}}_{{\bf \Sigma}}$, $i\in [N]$, $p\in [m_{vi}]$, $q\in [m_{zi}]$. Add an edge $(z_{iq},v_{ip})$ to ${\mathbb{T}}_{{\bf \Sigma}}$, then an input-unreachable $\lambda$-cycle emerges. According to Proposition \ref{proposition_3}, this indicates that the obtained NDS with the addition of edge $(z_{iq},v_{ip})$ is structurally uncontrollable. It further implies that the original NDS is structurally uncontrollable. Hence, together with Proposition \ref{proposition 4}, we have that the violation of either Condition (1) or Condition (2)  leads to structural uncontrollability of the NDS.  $\hfill\blacksquare$

{\begin{remark}  Corollary \ref{corollary_add} indicates that, the existence of an input-unreachable {{$\lambda$-edge}} implies the presence of either an FUM or a PDUM. This corollary is significant in proving the performance guarantee of the topology design procedure developed in the next section.
\end{remark}}

\section{Minimal Design of Subsystem Interconnection Topology}
In this section, as an application of the results established in Section V, we consider a minimal design of subsystem interconnection topology for an NDS to achieve structural controllability. All proofs of the results in this section are deferred to Appendix B.

\begin{spacing}{1.0}
Consider the following topology design problem.  Given $N$ linear time invariant subsystems ${\bf{\Sigma}}_1,...,{\bf{\Sigma}}_N$ whose subsystem dynamics are captured by (\ref{SubsystemDynamics}) with known first principle parameters $P^{(i)}|_{i=1}^N$, find the sparsest structured SCM $\Phi$, such that the constructed NDS is structurally controllable. Denote the NDS by ${\bf{{\bf \Sigma}}}(\Phi)$ when a structured $\Phi$ is designated.{\footnote{If ${\bf{{\bf \Sigma}}}(\Phi)$ is structurally controllable, a realization of $\Phi$ making the NDS controllable with probability $1$ can be obtained by setting values of the nonzero entries of $\Phi$ randomly from a set of real numbers with a sufficiently large cardinality. This is made clear by Proposition \ref{proposition 5} in Appendix A.}} Note that a nonzero entry of $\Phi$ is associated with an interconnection link between two subsystems, which will be called a subsystem link. For notational simplicity, the system matrices of the subsystem ${\bf \Sigma}_i$ with a fixed $P^{(i)}$, denoted by $A^{(i)}_{xx}, A^{(i)}_{xv}$, etc., are defined as  \end{spacing}  %\linespread{1.05} \selectfont{
{\footnotesize
\begin{equation}
\setlength{\abovedisplayskip}{2.0pt}
\left[\!\!\!{\begin{array}{*{20}{c}}
{A_{xx}^{(i)}}&{A_{xv}^{(i)}}&{B_{xu}^{(i)}}\\
{A_{zx}^{(i)}}&{A_{zv}^{(i)}}&{B_{zu}^{(i)}}\\
{C_{yx}^{(i)}}&{C_{yv}^{(i)}}&{D_{yu}^{(i)}}
\end{array}}\!\!\! \right] = {F_l}\left( {\begin{pmat}[{..|.}]
{A_{xx0}^{(i)}}&{A_{xv0}^{(i)}}&{B_{xu0}^{(i)}}&{E_1^{(i)}}\cr
{A_{zx0}^{(i)}}&{A_{zv0}^{(i)}}&{B_{zu0}^{(i)}}&{E_2^{(i)}}\cr
{C_{yx0}^{(i)}}&{C_{yv0}^{(i)}}&{D_{yu0}^{(i)}}&{E_3^{(i)}}\cr\-
{F_1^{(i)}}&{F_2^{(i)}}&{F_3^{(i)}}&{{H^{(i)}}}\cr
\end{pmat}},{P^{(i)}} \right),
\setlength{\belowdisplayskip}{3pt}
\end{equation}}and $m_{zi}$ and $m_{vi}$ denote the dimensions of internal output and internal input vectors of ${\bf \Sigma}_i$, respectively. All the remaining notations have the same definitions as those in the previous sections. Then, this problem can be formulated as follows.
\begin{problem}
 Given $N$ subsystems ${\bf{{\bf \Sigma}}}_1,...,{\bf{{\bf \Sigma}}}_N$,
\begin{equation}
{\centering \begin{array}{l}
\mathop {\min }\limits_{\Phi  \in \{0,*\}^{M_v\times M_z}} {\kern 1pt} {\kern 1pt} {\kern 1pt} {\kern 1pt} {\left\| \Phi  \right\|_0}\\
{\text{s.t.}}{\kern 1pt} {\kern 1pt} {\kern 1pt} {\kern 1pt} {\kern 1pt} {\kern 1pt} {\kern 1pt} {\kern 1pt} {\kern 1pt} {\kern 1pt}  {\bf{{\bf \Sigma}}} (\Phi ){\kern 1pt} {\kern 1pt} {\text{is structurally controllable}}
\end{array}}
\end{equation}
\label{prob2}
\end{problem}

A similar problem has been investigated in \cite{S_Pequito_2017}, \cite{Y_Zhang_2017}, for a structured system described by graphs. The above problem is different, noting that {{the dynamics of subsystems are numerically prescribed, and they can be high order and heterogeneous}}. Applications of this problem  include designing interaction topology for geographically distributed multi-agents to achieve consensus \cite{partitions}, designing communication links for geographically distributed sensors for data fusion.
\subsection{Feasibility and Complexity}
We first show the feasibility and complexity of Problem \ref{prob2}.

\begin{lemma} \label{Feasibility1}
Problem \ref{prob2} is feasible, if and only if the following three conditions are satisfied simultaneously: (i) for each $i\in [N]$,  $(A^{(i)}_{xx}, [B^{(i)}_{xu}\ A_{xv}^{(i)}])$ is controllable; (ii) $M_z\ge \max\nolimits_{1\le i \le m}{M_{ri}}$; (iii) there exists at least one integer $i\in [N]$ such that $G^{(i)}_{zu}(\lambda)\ne 0$, or for each $i\in [N]$, there is no $\lambda$-edge in the ACG ${\mathbb{T}}_i$.
\end{lemma}

Lemma {\ref{Feasibility1}} shows that to make an NDS structurally controllable, each of its subsystems should be controllable through its augmented input matrix $[B^{(i)}_{xu}\ A^{(i)}_{xv}]$. On the other hand, it also shows that even if an individual subsystem is not controllable through its external input matrix $B^{(i)}_{xu}$, the NDS could still be controllable through subsystem interconnections.   These
results are consistent with \cite{zhou_2015}, \cite{zhou_2018}, \cite{Y_Zhang_2016}.

\begin{proposition} \label{NP_hard1}
Problem \ref{prob2} is NP-hard.
\end{proposition}

\begin{algorithm}     \label{alg3}
\caption{{\small{Stage 1 of the topology design procedure: selecting minimal subsystem links to eliminate FUMs}}} % 算法的题目
\label{alg3} % 算法的标签
{{{\small{
\begin{algorithmic}[1]
\REQUIRE System matrices of subsystems ${\bf{{\bf \Sigma}}}_1$,...,${\bf{{\bf \Sigma}}}_N$.
\ENSURE   Approximated $\Phi$ for Problem \ref{prob3}

{\bf{Step 1: Use a greedy algorithm to approximate Problem \ref{prob4}}}

\STATE    Calculate $Y_i$, $Z_i$, for $i=1,...,{m}$;
 \STATE   Initialize $J \leftarrow \{1,...,M_v\}$, $J_{\rm grd}  \leftarrow \emptyset$;
 \WHILE{$g(J ) < \sum\nolimits_{i = 1}^{m} {{M_{ri}}} $}
\STATE ${s} \leftarrow a' \in \arg {\kern 1pt} {\kern 1pt} \max {\kern 1pt} {{\kern 1pt} _{a \in {J}\backslash J_{\rm grd} }}{\kern 1pt} {\kern 1pt} g(J_{\rm grd}  \cup \{ a\} ) - g(J_{\rm grd} )$;
\STATE $J_{\rm grd}  \leftarrow J_{\rm grd}  \cup \{{s}\}$;
\ENDWHILE

{\bf{Step 2: Use greedy coloring to construct a solution to Problem \ref{prob3} from $J_{\rm grd}$}}

\STATE  Find a collection $\{J_1,...,J_{m}\}$ of subsets of $[M_v+M_z]$ satisfying $|J_i|=M_{ri}$, $J_i\backslash {\cal B}\subseteq J_{\rm grd}$, such that ${\rm{rank}}(Q_{2i,J_i})=M_{ri}$, $|J_{i}\cap \cal{B}|$ is maximized per $i=1,...,m$,  where ${\cal B}=\{M_v+1,...,M_v+M_z\}$.\footnotemark[1]

\STATE Construct the coloring auxiliary graph ${\mathbb{G}}(J_1,...,J_{m})$;

\STATE  Initialize ${\mathbb{G}}_{\rm col}={\mathbb{G}}(J_1,...,J_{m})$, index $M_z$ colors as $1,...,M_z$, and use them to color ${\mathbb{G}}(J_1,..,J_{m})$ according to the following procedure:

 (i) for each vertex $j\in J_{\rm grd}\cap {\cal B}$, color vertex $j$ by the $(j-M_v)$-th color;

 (ii) after (i),  for each iteration, do the following operations, until there is no uncolored vertex in $\mathbb{G}_{\rm col}$:

\begin{itemize}
\item among all uncolored vertices, choose the one which is adjacent to the largest number of differently colored vertices, denoted by $v^\ast$;

\item if vertex $v^{\ast}$ has $M_z$ differently colored neighbors, assign $k^{\ast}_{\max}$ distinct colors to $v^{\ast}$, where \begin{equation} \label{color_k} k^{\ast}_{\max}=\max \{M_{ri}: v^\ast\in J_i, i=1,...,m\},\end{equation} and remove the edges between $v^{\ast}$ and its neighbors from $\mathbb{G}_{\rm {col}}$; otherwise, assign $v^\ast$ a color different from $v^\ast$'s colored neighbors, such that the number of already used colors is minimized;
\end{itemize}

\STATE  Map $\mathbb{G}_{\rm {col}}$ to $\Phi$: \quad $\Phi_{ij}=\ast$ if vertex $i$ is colored by color $j$ in $\mathbb{G}_{\rm {col}}$, $1\le i \le M_v$, $1 \le j \le M_z$; the rest entries of $\Phi$ are zero.
\end{algorithmic} }}}}
{\footnotemark[1]{\emph{{\footnotesize {This can be done efficiently in a simple greedy manner:  initialize $J_i \leftarrow \emptyset $;  in each iteration, choose an element $s$ from $\cal{B}$ and test ${\rm{rank}}(Q_{2i,J_i\cup\{s\}})-{\rm{rank}}(Q_{2i,J_i})$. If this value equals $1$ then $J_i \leftarrow J_i\cup\{s\}$; if ${\rm{rank}}(Q_{2i,J_i\cup\{s\}})-{\rm rank}(Q_{2i,J_i})$ keeps zero for all $s\in {{\cal{B}}\backslash {J_i}}$, choose $s\in {J_{\rm grd}}\backslash {\cal B}$ and continue the above iteration until ${\rm{rank}}(Q_{2i,J_i})=M_{ri}$.}}}}}
\end{algorithm}

\subsection{A Two-Stage Algorithm with Provable Approximation Bounds}

Since Problem \ref{prob2} is NP-hard, we propose a scalable algorithm to approximate it with some provable approximation bounds. This algorithm has two stages. Stage 1 is to select subsystem links to eliminate each FUM; Stage 2 is to add some extra subsystem links to guarantee the non-existence of a PDUM.

\subsubsection{Stage 1- FUM elimination}
In this Stage, starting from $N$ disconnected subsystems, we want to select the minimal number of subsystem links such that the obtained NDS ${\bf \Sigma}(\Phi)$ does not have an FUM. To clarify the dependence of $Q_1$ of Proposition \ref{proposition 4} on $\Phi$, rewrite $Q_1(\Phi) \triangleq [\Phi^{\intercal}\ I_{M_z}]$.  Define a function $f(\Phi): \{0,*\}^{M_v\times M_z}\rightarrow \mathbb{N}$ as %For two SCMs $\Phi_1, \Phi_2\in \{0,*\}^{M_v\times M_z}$, we say that $\Phi_1 \subseteq \Phi_2$, if $\Phi_{2ij}=0$ implies $\Phi_{1ij}=0$, $\forall i\in [M_v], j\in [M_z]$.
\begin{equation} f(\Phi ) = \sum\nolimits_{i = 1}^{m }{\rho ({\cal{M}}({Q_1}(\Phi )) \cap {\cal{M}}(Q_{2i}))} .\end{equation}
Suppose that the feasibility condition is satisfied. By Proposition \ref{proposition 4}, it suffices to see that the objective of the following Problem \ref{prob3} is equivalent to that of Stage 1.
\begin{problem}\label{prob3}
\begin{equation}
\begin{array}{l}
\mathop {\min }\limits_{\Phi   \in \{0,*\}^{M_v\times M_z} } {\kern 1pt} {\kern 1pt} {\kern 1pt} {\kern 1pt} {\kern 1pt} {\left\| \Phi  \right\|_0}\\
{\rm s.{\kern 1pt} {\kern 1pt} t.}{\kern 1pt} {\kern 1pt} {\kern 1pt} {\kern 1pt} {\kern 1pt} {\kern 1pt} {\kern 1pt} {\kern 1pt} {\kern 1pt} {\kern 1pt} {\kern 1pt} {\kern 1pt} {\kern 1pt} f(\Phi ) = \sum\nolimits_{i = 1}^{m} {{M_{ri}}}
\end{array}
\end{equation}
\end{problem}

\begin{lemma}\label{NP_hard2}
Problem \ref{prob3} is NP-hard.
\end{lemma}

Unfortunately, Problem \ref{prob3} is  NP-hard. Besides, it can be directly validated that $f(\Phi)$ is in general not submodular w.r.t. the nonzero elements of $\Phi$. The nonsubmodularity of $f(\Phi)$ might even prevent the existence of a nontrivial provable performance guarantee of using a simple greedy algorithm. Here, to seek for an algorithm {{with a provable performance bound}}, we propose the following alternative algorithm. This algorithm composes of two steps. The first step is to approximate a lower bound of  Problem \ref{prob3}, which can be formulated as a submodular function optimization problem.

To this end, let ${J}\subseteq [M_v]$ and denote by $\Omega({J}) \triangleq {\bf{{\bf {\bf diag}}}}\{I_{M_v,J}, I_{M_z}\}$. Define a function $g({J})$ as
\[g(J)= \sum\nolimits_{i = 1}^{m} {\rho ({\cal M}(\Omega (J)) \cap {\cal M}({Q_{2i}}))}\]
It can be directly proven that $g(J)=\sum\nolimits_{i = 1}^{m} {\rm rank}([{{Y_{iJ}}}\ {{Z_i}}])$. We then introduce another related problem, which is Problem \ref{prob4}, as follows.
\begin{problem} \label{prob4}
\[\begin{array}{l}
\mathop {\min }\limits_{J \subseteq [{M_v}]} {\kern 1pt} {\kern 1pt} {\kern 1pt} {\kern 1pt} {\kern 1pt} |J|\\
{\rm s.{\kern 1pt} {\kern 1pt} t.{\kern 1pt} {\kern 1pt}} {\kern 1pt} {\kern 1pt} {\kern 1pt} {\kern 1pt} {\kern 1pt} {\kern 1pt} {\kern 1pt} {\kern 1pt} {\kern 1pt} {\kern 1pt} {\kern 1pt} g(J) = \sum\nolimits_{i = 1}^{m} {{M_{ri}}}
\end{array}\]
\end{problem}

\begin{proposition} \label{Submodular_proposition}
$g({J})$ is submodular on ${J} \subseteq [M_v]$. Moreover, denote the solution of Problem \ref{prob4} returned by the greedy algorithm stated as Step 1 in Algorithm {\ref{alg3}} by ${J}_{{\rm grd}}$, and the optimal solution to Problem \ref{prob3} by $\Phi^*$. Then, it holds that {\small{
\[\frac{{|{J_{\rm grd}}|}}{{{{\left\| {{\Phi ^*}} \right\|}_0}}} \le 1 + \log \frac{M_{\rm def}}{{\sum\nolimits_{i = 1}^{m} {{M_{ri}} - g({J_{T - 1}})} }} \le 1 + \log {M_x},\]}}in which
$M_{\rm def}\triangleq {\sum\nolimits_{i = 1}^{ m} {\sum\nolimits_{j = 1}^{N} (m_{x{j}}- {\rm rank}([\lambda_i I- A_{xx}^{(j)}\ B_{xu}^{(j)}]))}}$, ${J}_{T-1}$ is the return value of the second-to-last iteration of the greedy algorithm.
\end{proposition}

From its definition, $M_{\rm def}$ is obviously the total row rank deficiency of all subsystems at their uncontrollable modes.

The second step uses the greedy coloring techniques \cite{DB_West_graph} to restore a feasible solution to Problem \ref{prob3} from $J_{\rm grd}$. Graph coloring is the problem of coloring vertices of a graph such that any two adjacent vertices do not share the same color. Greedy coloring is a heuristic method towards this problem that assigns to a vertex with {{the smallest available color among all colored vertices not used by its neighbours in a specific order, adding a new color if needed}}. The adoption of graph coloring is motivated by the following observations. Notice that for a given $\Phi$,  ${\rho ({\cal{M}}({Q_1}(\Phi)) \cap {\cal{M}}(Q_{2i}))}=M_{ri}$ implies that there are $M_{ri}$ nonzero entries located in the columns of $Q_1(\Phi)$ indexed by $J_i$ for some $J_i \subseteq [M_z+M_v]$ with $|J_i|=M_{ri}$ and ${\rm rank}(Q_{2i,J_i}) = M_{ri}$, such that any two of them are in different rows of $Q_1(\Phi)$.
Regrading the column indices of these nonzero entries as vertices and their row indices as colors, an equivalent statement of the aforementioned condition on the SCM $\Phi$ is that, the clique formed by these $M_{ri}$ vertices are colored such that any adjacent vertices do not share a common color, i.e., an $M_{ri}$-coloring \cite{DB_West_graph}. It should be noted that a similar idea is also adopted in \cite{Y_Zhang_2018} for a different yet simpler problem.

To formulate the greedy coloring process, for a given collection $\{J_1,...,J_{m}\}$ of subsets of $[M_v+M_z]$ with $|J_i|=M_{ri}$, define the coloring auxiliary graph as ${\mathbb{G}}(J_1,...,J_{m})=({\cal{V}}_{\rm{{col}}},{\cal{E}}_{\rm{{col}}})$, with ${\cal{V}}_{\rm{{col}}}=\bigcup\nolimits_{i = 1}^{m} {{J_i}}$ and ${\cal{E}}_{\rm{{col}}}=\left\{ {({j_1},{j_2}):{j_1},{j_2} \in {J_i},{j_1} \ne {j_2}, i=1,...,m} \right\}$. Equivalently, ${\mathbb{G}}(J_1,...,J_{m})$ is the union of $m$ individual cliques formed by vertices of $J_i$, $i=1,...,m$. A precise statement for this procedure is given in Algorithm {\ref{alg3}}. The greedy coloring procedure is in Substep $9$ of Algorithm \ref{alg3}, where we adopt a dynamically updated ordering, and as stated in (\ref{color_k}), permit one vertex to be assigned with more than one color. A bound for its performance is given in Theorem \ref{theorem_last}.  An illustrative example is given in Section VII.

\subsubsection{Stage 2-PDUM elimination} In this stage, the objective is to add minimal subsystem links to eliminate all PDUMs of the NDS obtained after Stage 1.  The associated procedure is stated as Algorithm \ref{alg1}. In this procedure,  an SCC is called input-unreachable source SCC, if in the associated n-ACG, this SCC is input-unreachable, and there is no incoming edge from any other SCCs to any of its vertices. Let $p_{\rm {\rm {ius}}}$ be the number of input-unreachable source SCCs that contain a $\lambda$-edge in the n-ACG obtained after Stage 1, denoted by ${\mathbb{T}}_{{\bf{{\bf \Sigma}}}{\text{s1}}}$. Then, it can be seen that, at least $p_{\rm {\rm {ius}}}$ subsystem
links are needed to guarantee the input-reachability of all $\lambda$-cycles in ${\mathbb{T}}_{{\bf{{\bf \Sigma}}}{\text{s1}}}$. Meanwhile, Algorithm \ref{alg1} provides a
constructive procedure on how $p_{\rm {\rm {ius}}}$ subsystem links are sufficient to do so. From Proposition \ref{proposition_3}, $p_{\rm {\rm {ius}}}$ is the minimal number of subsystem links whose addition to ${\mathbb{T}}_{{\bf{{\bf \Sigma}}}{\text{s1}}}$ eliminates each PDUM of the NDS obtained after Stage 1.

The following theorem gives some approximation bounds for the two-stage topology design procedure.

{\begin{theorem} \label{theorem_last}
 The two-stage topology design procedure overall returns an $O(2M_{\rm rmax}{\rm log}(M_{\rm def}))$ approximation for Problem \ref{prob2}, where $M_{\rm rmax}\triangleq \max \nolimits_{1\le i \le m} M_{ri}$. In Stage 1, Algorithm \ref{alg3} returns an $O(M_{\rm rmax}{\rm log}(M_{\rm def}))$ approximation for Problem \ref{prob3}. In addition, if every vertex of $\mathbb{G}_{{\rm col}}$ in Substep $10$ of Algorithm \ref{alg3} is colored by only one color, then the approximation factor becomes $O({\rm log}(M_{\rm def}))$ for Problem \ref{prob3}.
\end{theorem}}

There are two promising features in the above topology design procedure. First, {{it has low computational complexity.}} To be specific, Step 1 of Stage 1 has complexity $O(NM_{\rm def})$, which is no more than $O(M^2_x)$; particularly, the call of a matroid intersection subroutine is no longer needed. Step 2 of Stage 1 incurs at most an $O(M_z^3)$ complexity. Stage 2 has an $O(M_vM_z)$ complexity. Hence, the overall time complexity is $O(M_z^3+M^2_x)$. Note that we usually have $M_z\ll M_x$ in actual applications \cite{zhou_2015}. This means that the time complexity of our approach usually increases quadratically with $M_x$ in practice. Second, {{it has a provable approximation bound}}, while other possible heuristics might not. The approximation bound $2M_{\rm rmax}{\rm log}(M_{\rm def})$ might seem loose. However, when the subsystems have smaller dimensions of uncontrollable subspaces which leads to a smaller $M_{\rm def}$, or have more heterogeneous eigenvalues which leads to a smaller $M_{\rm rmax}$, this bound becomes tighter. Moreover, if $M_z\gg M_{\rm rmax}$, the approximation bound becomes $2{\rm log}(M_{\rm def})$, since in this case every vertex of ${\mathbb{G}}_{\rm col}$ will be colored by only one color. In fact, the proof of Proposition \ref{NP_hard1} indicates that approximating Problem \ref{prob2} within a factor $c{\rm log}(M_{\rm def})$ is NP-hard for any constant $c<1$.  The numerical example in Section VII shows that our approach can sometimes return the optimal solution.

The above topology design procedure may also be applied in the following situation. Given a collection of disconnected subsystems, adding the minimal subsystem links to construct an NDS without unstable FUMs, i.e., the FUM with a nonnegative real part. This objective is closely related to constructing an NDS that can be stabilized by state feedback \cite{Kailath_1980}.\footnote{If there is no input-unreachable $\lambda$-cycle in the n-ACG of the obtained NDS, the NDS will have no PDUM. In this case, a sufficient condition for stabilization can be obtained.} To this end, only the unstable FUMs are needed to be considered. Hence, a slightly modified version of Algorithm \ref{alg3} may be suitable for this problem. The numerical example in Section VII illustrates this capability.

%(ii) Given a collection of subsystems, possibly including unstable subsystems, design the minimal subsystem interconnection topology to construct an NDS without
%unstable fixed modes. This is closely related to constructing a stable NDS using possibly unstable subsystems. For this purpose,  set $B^{(i)}_{zu}\equiv 0$, $B^{(i)}_{xu}\equiv 0$ per $i\in [N]$. Then, a modified version of Algorithm \ref{alg3} may also be suitable to approximate it.
{{\small{
\begin{algorithm}
\caption{{{\small {Stage 2 of the topology design procedure: adding minimal subsystem links to eliminate PDUMs}}}} % 算法的题目
\label{alg1}
\begin{algorithmic}[1]
{{{\small
\REQUIRE The associated n-ACG after Stage 1, denoted by ${\mathbb{T}}_{{\bf{{\bf \Sigma}}}{\text{s1}}}$.

\ENSURE   Approximated $\Phi$ for Problem \ref{prob2}
\STATE   Identify the input-unreachable source SCCs of ${\mathbb{T}}_{{\bf{{\bf \Sigma}}}{\text{s1}}}$,  and among them denote those that have a $\lambda$-edge by ${\hat {\cal N}}_1,..., {\hat {\cal N}}_{p_{\rm {\rm {ius}}}}$;  ${\mathbb{T}}_{{\bf \Sigma}} \leftarrow {\mathbb{T}}_{{\bf{{\bf \Sigma}}}{\text{s1}}}$;
\FOR{$i=1,...,p_{\rm {\rm {ius}}}$}
\STATE Update ${\mathbb{T}}_{{\bf \Sigma}}$ by adding a link from $z_l$ to $v_j$, where the internal input vertex $v_j$ belongs to $\hat {\cal N}_i$, the internal output vertex $z_l$ belongs to the input-reachable SCCs of ${\mathbb{T}}_{{\bf \Sigma}}$ {{updated}} in the previous step.
\ENDFOR
\STATE Map ${\mathbb{T}}_{{\bf \Sigma}}$ to $\Phi$ according to the correspondence between  ${\mathbb{T}}_{{\bf \Sigma}}$ and $\Phi$. }}}
\end{algorithmic}
\end{algorithm}
}}}
\vspace{0.1cm}
\section{An Illustrative  Example}

In this section, a numerical example is given to illustrate the results obtained in the previous sections.

\begin{figure}
  \centering
  \includegraphics[width=2.75in]{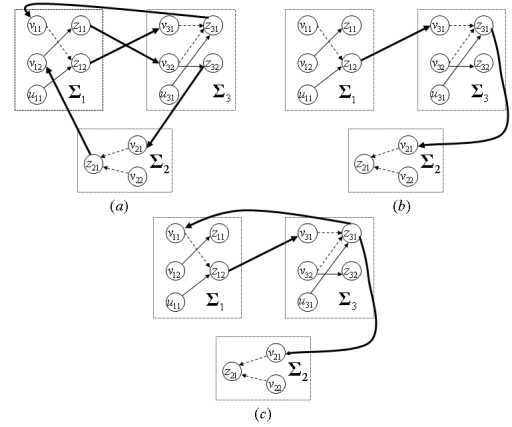}
  \caption{Illustration of n-ACGs of the example in Section VII. The dotted edges represent $\lambda$-edges, while the bold ones represent subsystem links.}\label{n-ACG}
\end{figure}

\begin{figure}
  \centering
  \includegraphics[width=3.6in]{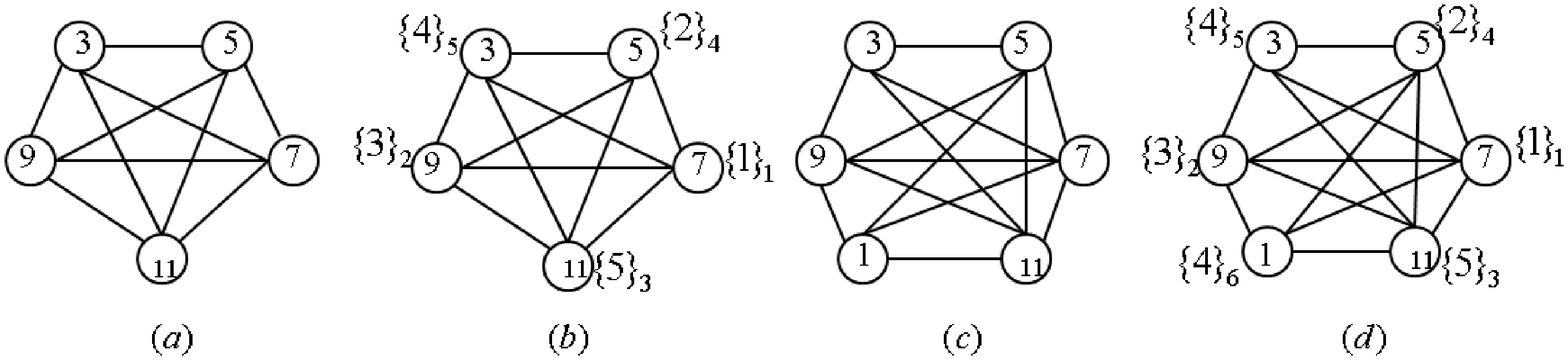}
  \caption{The $\mathbb{G}(J_1,...,J_l)$ ((a) and (c) are respectively for $l=2$ and $3$) before coloring and $\mathbb{G}_{{\rm col}}$ ((b) and (d)) after coloring. The number in the bracket near each vertex is the color assigned to this vertex, and the subscript number of the bracket is the order that this vertex is colored in Algorithm \ref{alg3}.}\label{coloring_process}
\end{figure}

Consider an NDS consisting of $3$ subsystems. The system matrices for these subsystems are respectively
{\footnotesize{
$Pa({{\bf \Sigma} _1}) = \begin{pmat}[{.||.}]
0&0& 1& 1&0\cr
0&{-1}& 0& 1&0\cr\-
0&0& 0& 0&1\cr
1&1& 0& 0&0\cr
\end{pmat}$, $Pa({{\bf \Sigma} _2}) = \begin{pmat}[{.||.}]
1&0& 0& 1&0\cr
2&0& 0& 1&0\cr\-
1&0& 0& 1&0\cr
\end{pmat}$, $Pa({{\bf \Sigma} _3}) = \begin{pmat}[{.||.}]
{ - 1}&0& 0& 1&1\cr
2&1& 0& 0&0\cr\-
0&1& 1& 1&1\cr
0&0& 0& 0&1\cr
\end{pmat},$}}
where {\footnotesize{ $Pa({{\bf \Sigma} _i}) \triangleq \begin{pmat}[{||}]
{A_{xx}^{(i)}}& {B_{xu}^{(i)}}& {A_{xv}^{(i)}}\cr\-
{A_{zx}^{(i)}}& {B_{zu}^{(i)}}& {A_{zv}^{(i)}}\cr
\end{pmat}$}}. Let a structured SCM $\Phi_{a} =${{$\tiny{\begin{pmat}[{.||.}]
{0}&{0}& {0}& *&{0}\cr
{0}&{0}& *& {0}&{0}\cr\-
{0}&{0}& {0}& {0}&*\cr
{0}&{0}& {0}& {0}&{0}\cr\-
{0}&*& {0}& {0}&{0}\cr
*&{0}& {0}& {0}&{0}\cr
\end{pmat}}$}}.

First, consider the SCM $\Phi_{a}$.  The corresponding n-ACG $\mathbb{T}_{{\bf \Sigma}}$  is given in Fig. \ref{n-ACG}(a). From  Fig. \ref{n-ACG}(a), there is an input-unreachable $\lambda$-cycle in $\mathbb{T}_{{\bf \Sigma}}$, which is $\{v_{12}\to z_{11}\to v_{32} \to z_{32}\to v_{21} \to z_{21} \to v_{12}\}$. Hence, by Proposition \ref{proposition_3}, the above NDS has at least one PDUM, which means that it is not controllable under the above structured SCM {{$\Phi_a$}} with any specific link weights. This can be validated by some algebraic manipulations.

Second, we show the application of the modified version of Algorithm 1 in selecting the minimal subsystem links to construct an NDS without unstable FUMs.  Note that the matrix pair $(A^{(i)}_{xx},B^{(i)}_{xu})$ has some  unstable uncontrollable modes for $i=2,3$. This means that Subsystems ${\bf \Sigma}_2$ and ${\bf \Sigma}_3$ cannot be stabilized by state feedback when isolated.  The set of subsystem eigenvalues is ${\bf \Lambda}=\{1,0,-1\}$, with unstable mode set being $\hat {\bf \Lambda}= \{1, 0\}$. Let $\lambda_1=1$, $\lambda_2=0$. Applying Algorithm \ref{alg3} by letting $m=2$, after Step 1, we get $J_{\rm grd}=\{3,5\}$. For Step 2 in Substep 7 of Algorithm \ref{alg3}, we have $J_1=\{3,5,7,11\}$ and $J_2=\{3,5,7,9,11\}$. The coloring auxiliary graph ${\mathbb{G}}({{ J_1,J_2}})$ is given in Fig. \ref{coloring_process}(a), with its associated colored graph ${\mathbb {G}}_{\rm {col}}$ being Fig. \ref{coloring_process}(b) obtained by executing Step 2 of Algorithm \ref{alg3}. According to the correspondence between coloring and the SCM, a structured $\Phi$ with its $(3,4)$-th and $(5,2)$-th entries being nonzero is obtained. That is, two subsystem links $(z_{31},v_{21})$ and $(z_{12},v_{31})$, as illustrated in Fig. \ref{n-ACG}(b), are sufficient to eliminate these unstable FUMs. From Fig. \ref{n-ACG}(b), there does not exist an input-unreachable $\lambda$-cycle. Hence, there is no PDUM in the NDS associated with Fig. \ref{n-ACG}(b). Since any other FUM if there exists, is stable, for ``almost each'' numerical $\Phi_0$ with the structure shown in Fig. \ref{n-ACG}(b), the associated NDS ${\bf \Sigma}(\Phi_0)$ is stabilizable by a state feedback. To validate it, assign a random weight to each nonzero entry of the SCM $\Phi$. Using the PBH test, it is observed that every uncontrollable mode of the obtained NDS is fixed at $\lambda=-1$, which means that the NDS can be stabilized.

Now suppose that our goal is to find the sparsest SCM to construct a structurally controllable NDS. Continuing the above procedure, let $\lambda_3=-1$.
Following a similar procedure described above, we get $J_{\rm grd}=\{1,3,5\}$ in Step 1 of Algorithm \ref{alg3} with $J_1=\{3,5,7,11\}$, $J_2=\{3,5,7,9,11\}$, $J_3=\{1,5,7,9,11\}$ in Substep 7 of Algorithm \ref{alg3}. The associated graphs ${\mathbb{G}}({J_1,J_2,J_3})$ and ${\mathbb G}_{\rm col}$ before and after coloring are given respectively in Fig. \ref{coloring_process}(c) and Fig. \ref{coloring_process}(d). The corresponding interconnection topology is illustrated in Fig. \ref{n-ACG}(c), which shows that there does not exist an input-unreachable $\lambda$-cycle. Hence, the obtained SCM with its $(5,2)$-th, $(1,4)$-the and $(3,4)$-th entries being nonzero, makes the NDS structurally controllable. Through an exhaustive search, it becomes clear that this interconnection topology is optimal in the sense that it has the minimal number of subsystem links making the associated NDS structurally controllable. %, which are the $(5,2)$, $(1,4)$ and $(3,4)$ positions of $\Phi$

\section{Conclusions}
This paper investigates structural controllability of an NDS, in which unknown parameters are allowed to exist in both subsystem dynamics in an LFT parameterized way and in subsystem interconnections, and each subsystem may have high-order, heterogeneous dynamics. Some results are first obtained about structural controllability of an LFT parameterized plant under a diagonalization assumption, which further lead to some necessary and sufficient conditions respectively for the NDS to have an FUM, a PDUM, and to be structurally controllable. These conditions can be verified efficiently, and give some intuitive insights on how the network controllability is influenced by subsystem input-output relations, subsystem uncontrollable modes and subsystem interconnection topology. Based on them, a minimal design of subsystem interconnection topology is considered for an NDS to achieve structural controllability. A two-stage algorithm is proposed with some provable performance bounds. Further researches include studying more generic properties for an NDS with unknown parameters, like the generic dimension of controllable subspaces, structural controllability with a random switching topology, etc.

\begin{appendices}
\section{General complexity and hardness of structural
controllability verification: high-rank case}
In this appendix, we discuss the general computational complexity and hardness of structural controllability verification when the coefficient
matrix of a variable in system matrices is not restricted to be rank-one.
Specifically, for the pair in (\ref{LFT parameterized}) we set $A_{zv}\equiv 0$ and equivalently write the following $(A,B)$ as in \cite{Anderson_1982}:
\begin{equation} \label{linear A_B} A = {A_0} + \sum\nolimits_{i = 1}^k {{s_i}{A_i}} ,B = {B_{0}} + \sum\nolimits_{i = 1}^k {{s_i}{B_i}},\end{equation} where $A_i\in\mathbb{R}^{n\times n}$ and $B_i\in \mathbb{R}^{n\times q}$ are constant matrices, and $s_1, ..., s_k$ are free parameters.

\begin{definition}[\rm RP, \cite{computation_complexity_modern_2009}] \label{definition rp}
Randomized polynomial time (RP) is the complexity class of problems for which a random algorithm exists with:
(i) it always runs in polynomial time in the input size;
(ii) if the correct answer is YES, it returns YES with probability at least $1/2$; if the correct answer is NO,  then  it always returns NO.
\end{definition}

\begin{lemma}[\cite{Separation_1980}] \label{lemma_A1}
Let $f(s_1,...,s_t)$ be a multivariate polynomial of variables $s_i|_{i=1}^t$ with real coefficients that is not identically zero. Let $\mathbb{V}$ be an arbitrary subset of $\mathbb{R}$ with a finite cardinality, $d$ the total degree of $f(s_1,...,s_t)$ (i.e., the highest degree of one monomial), and $\mathbb{V}^{t}$ the set of all $t$-element arrays with each element chosen from $\mathbb{V}$. Then, if $d\le |\mathbb{V}|$, $f(a)\ne 0$ for at least $(|\mathbb{V}| - d)|\mathbb{V}|^{t - 1}$ elements out of all $a\in \mathbb{V}^{t}$.
\end{lemma}

From Lemma \ref{lemma_A1}, the following proposition is obtained.

\begin{proposition} \label{proposition 5} {Verifying the structural controllability of $(A,B)$ in (\ref{linear A_B}) is {\rm {RP}}}. More specifically, let ${r_{\max }}\triangleq \mathop {\max {\kern 1pt} }\nolimits_{1 \le i \le k} {\kern 1pt} {\kern 1pt}{\rm rank}([{A_i}\ {B_i}])$, and the set of real numbers $\mathbb{V}\subseteq \mathbb{R}$ satisfies $|\mathbb{V}|= \min \{ 2kn{r_{\max}}, 2n^2 \}$.  If $(A,B)$ is structurally controllable, then, randomly choose an element $(s_1,...,s_k)\in \mathbb{V}^{k}$, with probability at least $1/2$, the obtained numerical $(A,B)$ of (\ref{linear A_B}) is controllable.
\end{proposition}

\emph{Proof:}  A numerical $(A,B)$ is controllable if and only if the matrix $\mathcal{C}$ defined in (\ref{Rocenber}) with dimension ${n^2} \times n(n + q - 1)$ is of FRR \cite{Kailath_1980}.
{\footnotesize{\begin{equation} \label{Rocenber} {\cal C}\triangleq \left[ \begin{matrix}
B&I&0&0&0&0&0\\
0&{ - A}&B&I&0&0&0\\
 \vdots &{}& \ddots & \ddots & \ddots &{}& \vdots \\
0&0&0&{ - A}&B&I&0\\
0&0&0&0&0&{ - A}&B
\end{matrix} \right]\end{equation}}}If $(A,B)$ in (\ref{linear A_B}) is structurally controllable, there exists at least $n^2$ columns of $\mathcal{C}$, the set of whose column indexes is denoted by $\mathcal{K}$, such that ${\rm det}\, {\cal C}_{\mathcal{K}}$ cannot be identically zero. Let $d$ be the total degree of ${\rm det}\,{\cal C}_{\mathcal{K}}$. It's obvious that $d\le n^2$. Notice that also, the degree of $s_i$ in ${\rm det}\,{\cal C}_{\mathcal{K}}$ is bounded by
${\rm{rank}}({B_i}) + (n - 1){\rm{rank}}([{A_i}\ {B_i}]) \le n{r_{\max }}$, thus $d \le knr_{\max}$. According to Lemma \ref{lemma_A1}, choose randomly an element $(s_1,...,s_k)$ from $\mathbb{V}^{k}$, then under the condition $|\mathbb{\mathbb{V}}| \ge 2d$, the probability of ${\rm det}\,{\cal C}_{\mathcal{K}}\ne 0$, denoted by  $\bf{Pr}$, satisfies

\begin{equation}
\label{probability}
{\bf{Pr}} \ge {\frac{{\left( {|\mathbb{V}| - d} \right)|\mathbb{V}{|^{k - 1}}}}{{|\mathbb{V}{|^k}}}} \ge \frac{1}{2}.
\end{equation}
Note that verifying whether $(A,B)$ in (\ref{linear A_B}) for a given $(s_1,...,s_k)$ is controllable can be done in polynomial time by testing the rank of the associated controllability matrix. By Definition \ref{definition rp}, the result follows.  $\hfill\blacksquare$

From Proposition \ref{proposition 5}, it can be seen that, if the cardinality of $\mathbb{V}$ tends to be sufficiently large, the probability in (\ref{probability}) will tend to be $1$. In other words,  Proposition \ref{proposition 5} quantifies the statement that controllability is a generic property by measuring the cardinality of the parameter space and the probability of controllability for a randomly chosen system. On the other hand, under the generally believed conjectures \cite{computation_complexity_modern_2009} that $\rm RP {\rm =} \rm P$ and ${\rm P} \ne {\rm NP}$, it indicates that {verifying the structural controllability of $(A,B)$ of (\ref{linear A_B}) should not be NP-hard}, and there should be some derandomized algorithm that can do this efficiently. However, the following theorem reveals the hardness of finding such a deterministic algorithm. In this theorem, the definitions of complexity classes $\rm {NEXP}$ and $\rm P\setminus\rm poly$, could be referred to e.g.,  \cite{computation_complexity_modern_2009} and \cite{A_Shpilka_survey}, and $\rm {Permanent}$ of an $n\times n$ matrix $M$ is the sum over all permutations of $n$ columns, of a product of $n$ terms, where the $i$-th term in the product is the term in the $i$-th row and the permutation of the $i$-th column.

\begin{theorem} \label{theorem3}
If one can deterministically verify whether the system in (\ref{linear A_B}) is structural controllable in polynomial time (or even in sub-exponential time), then either  (i)  $\rm {NEXP} \not\subset \rm P\setminus \rm poly $ or (ii) $\rm {Permanent}$ is not computable by polynomial sized arithmetic circuits over $\mathbb{Z}$.
\end{theorem}

Theorem \ref{theorem3} means that finding a deterministic algorithm to verify structural controllability of (\ref{linear A_B}) in polynomial time or even in sub-exponential time, is at least as hard as proving
Statement (i) or (ii) in this theorem, which are two open problems in arithmetic circuits \cite{A_Shpilka_survey}. Although it is generally believed that $\rm {NEXP} \not\subset \rm P\setminus \rm poly $ and $\rm {Permanent}$ requires super-polynomial sized circuits \cite{computation_complexity_modern_2009}, it is also commonly agreed that we are far away from proving these lower bounds \cite{A_Shpilka_survey}.
Nevertheless, for the case where $[A_i\ B_i]$ has an arbitrary rank, some efficient random algorithms, or black-box polynomial identity testing algorithms \cite{A_Shpilka_survey, Separation_1980}, could be adopted.

To prove Theorem \ref{theorem3}, we leverage a special problem about arithmetic circuit complexity, called the {{symbolic determinant identity testing problem (SDIT)}} \cite{Derandom_2004}. This problem can be stated as follows. For a square matrix $M$ with its entries being either a constant integer or a variable, where each variable could appear more than once, determine whether or not its determinant is identically zero.  By reducing the general SDIT to a special instance of structural controllability verification problem, we give the proof with the following Lemma \ref{lemma_A2}. Notice that Theorem \ref{theorem3} is not contradictory to the main results of this paper, as $\rm {SDIT}$ for rank-one case is deterministically solvable in polynomial time \cite{A_Shpilka_survey}, while is more challenging and still open for the general high-rank case.
\begin{lemma}[\cite{Derandom_2004}, \cite{A_Shpilka_survey}] \label{lemma_A2}
The following three statements cannot be simultaneously true: (i) $\rm {SDIT}$ can be solvable in polynomial time (or even in sub-exponential time); (ii)  $\rm {NEXP} \subset \rm P\setminus \rm poly $;  (iii) $\rm {Permanent}$ is computable by polynomial sized arithmetic circuits over $\mathbb{Z}$.
\end{lemma}

{\emph{Proof of Theorem \ref{theorem3}:}} Let $B$ be an $n\times n$ matrix, whose entries is either a constant integer or a variable. Let $A=I_{n}$. It's obvious that $(A,B)$ has the form of (\ref{linear A_B}), and all the coefficient matrices can be obtained in polynomial time. By the PBH test, $(A,B)$ is structurally controllable, if and only if $\det B$ is not identically zero. This means that verifying the structural controllability of $(A,B)$, is at least as hard as the SDIT on $B$. The result then follows from Lemma \ref{lemma_A2}. $\hfill\blacksquare$

\section{Proofs of some technical results}

{\emph{Proof of Lemma {\ref{lemma_proof_mid_2}}:}} Without loss of generality, denote cycles $C_1$ and $C_2$ by ${C_1} = \left\{ {{i_1} \to {i_2} \to  \cdots  \to {i_p} \to {i_1}} \right\}$, ${C_2} = \left\{{{i_1} \to i_2' \to  \cdots  \to i_p' \to {i_1}}\right\}$, where $\{i_1,i_2,\cdots,i_p\}$ and $\{i_1,i'_2,\cdots,i'_p\}$ are two distinct permutations of $\{1,\cdots,p\}$.
Let $t=\min\{j\in \{2,\cdots,p\}: i_{j}\ne i'_{j}\}$. Then, there must exist some integer $t'$ such that $i'_{t'}=i_t$ while $t< t'\le p$. Notice that $(i_{t-1},i'_{t'})\in {E}(C_1)$. Hence, the cycle ${C^*} = \left\{ {{i_1} \to {i_2} \to  \cdots  \to {i_{t - 1}} \to i_{{t'}}' \to i_{{t'} + 1}' \to  \cdots  \to i_p' \to {i_1}} \right\}$ contains edge $(i_1,i_2)$ and has a length no more than $p-1$. Similar analysis on edge $(i_1,i_2)$ can be applied to any other edge of $C_1$. $\hfill\blacksquare$

{{\emph{Proof of Lemma {\ref{lemma_proof_2}}:}}} For the ease of notation and without loss of generality, let us consider ${\cal J}=\{1,...,n_s\}$ and denote ${S}_{\cal J}=\{s_1,...,s_{n_s}\}$. Define a {weighted multigraph} associated with the TFM $G^\mathcal{J}_{zv}(\lambda)P_{\mathcal{J}}-I$ as ${\mathbb D}^{\mathcal{J}}=({\cal V}_{\mathcal{J}},{\cal E}_{\mathcal{J}},W_{\mathcal{J}})$, where the vertex set ${\cal V}_{\mathcal{J}}=\{1,...,n_s\}$, the edge set ${\cal E}_{\mathcal{J}}={\cal E}^0_{\mathcal{J}}\bigcup {\cal E}_{I}$ in which ${\cal E}^0_{\mathcal{J}} = \{(j,i):[{G^{\mathcal{J}}_{zv}}(\lambda )]_{ij} \ne 0\}$ and ${\cal E}_{I}=\{(i,i):i=1,...,n_s\}$, and the weight set {\small{$W_{\mathcal{J}} = \{{w(e):w(e) = [{G^{\mathcal{J}}_{zv}}{{(\lambda )}]_{ij}}{s_j}\, {\rm for}\, e=(j,i) \in {\cal E}^0_{\mathcal{J}}, w(e)=-1\, {\rm for}\, e \in {\cal E}_{I}}\}$}}. Moreover, for an edge $e=(i,j) \in {\cal E}^0_{\mathcal{J}}$, define ${\bar w}(e)= [{G^{\mathcal{J}}_{zv}}{(\lambda )}]_{ij}$. From this definition, some vertices may have multiple self-loops, which leads to the term ``multigraph''. Let $\tilde S_{n_s}$ be the collection of all matchings of ${\mathbb D}^{\mathcal{J}}$ with size $n_s$, and denote the $i$-th matching by $C_i$. It is well-known that each matching $C_i$ corresponds to a collection of disjoint cycles that span ${\mathbb D}^{\cal J}$ \cite{DB_West_graph}.  From the definition of determinant \cite{Matrix analysis}, we have
\begin{equation}
\label{determinant}
{\rm{det}} ({G^{\mathcal{J}}_{zv}}{(\lambda )}{P_{\mathcal{J}}} - I) = \sum\nolimits_{{C_i} \in  \tilde S_{n_s}} {{\mathop{\rm sgn}} ({C_i})} \prod\nolimits_{e \in
{C_i}} {w(e)},
\end{equation}
where ${\rm sgn}({C_i}) \in \{-1,1\}$ is the sign associated with  $C_i$.

Notice that a nonzero term $\prod\nolimits_{e \in {C_i}}{w(e)}$ in (\ref{determinant}) could either be {a monomial of $\{s_1,...,s_{n_s}\}$ depending on $\lambda$, a monomial of $\{s_1,...,s_{n_s}\}$ independent of $\lambda$, or the constant $(-1)^{n_s}$}. The {{necessity}} is then obvious, since if there dose not exist a $\lambda$-cycle in $\mathbb{L}^{\mathcal{J}}$ (therefore neither in ${\mathbb D}^{\mathcal{J}}$), there exists no term in ${\rm{det}} ({G^{\mathcal{J}}_{zv}}{(\lambda )}{P_{\mathcal{J}}} - I)$ that depends on the variable $\lambda$.

 To show the {{sufficiency}}, notice that in the digraph ${\mathbb D}^{\mathcal{J}}$, all outgoing edges from a common vertex $i \in {\cal V}_{{\cal J}}$ have weights with the same factor $s_i$ (excluding the self-loop with weight $-1$). This leads to the fact that for two distinct matchings $C_i$ and $C_j$, their associated terms $\prod\nolimits_{e \in {C_i}}{w(e)}$ and $\prod\nolimits_{e \in  {C_j}}{w(e)}$ could cancel each other out. Suppose that there exists at least one  $\lambda$-cycle in $\mathbb{L}^{\mathcal{J}}$, and so does in ${\mathbb D}^{\mathcal{J}}$ by definition. Denote the $\lambda$-cycle which has the {{shortest length}} among  all these $\lambda$-cycles in ${\mathbb D}^{\mathcal{J}}$ by ${\bar C_{\min}}=\{ {i_1} \to \cdots \to {i_{{k_{\min }}}} \to {i_1}\}$, where $k_{\min}$ is the length of $\bar C_{\min}$. Notice that ${\bar C_{\min}}$ may not be unique. Denote by $\tilde C_{\min}$ the matching in ${\mathbb D}^{\mathcal{J}}$ constituted by $\bar C_{\min}$ and the rest $n_s-k_{\min}$ self-loops with weights $-1$. We declare that the term associated with $\tilde C_{\min}$, expressed in (\ref{minimum_term}), where $k_{\min}+1$ is defined to be $1$ and the sign is ignored, cannot be varnished by other terms in ${\rm{det}} ({G^{\mathcal{J}}_{zv}}{(\lambda )}{P_{\mathcal{J}}} - I)$.
\begin{equation}\label{minimum_term}  \prod\nolimits_{j = 1}^{{k_{\min }}} {\bar w(({{i_j},{i_{j + 1}}}))} \prod\nolimits_{j = 1}^{{k_{\min }}} {{s_{{i_j}}}}\end{equation}

To show this, first observe that there does not exist a set of vertices ${{\cal Z}}_s^c \subseteq {\rm{\{ }}1,...,{n_s}{\rm{\} \backslash }}\{ {i_1},...,{i_{{k_{\min }}}}{\rm{\} }}$, such that ${{\cal Z}}_s^c \cup \{ {i_1},...,{i_{{k_{\min }}}}\}$ forms a larger cycle containing $\bar C_{\min}$ and corresponds to a monomial of the form $f(\lambda)\prod\nolimits_{j = 1}^{{k_{\min }}} {{s_{{{i}_j}}}}$, where $f(\lambda) \in F(\lambda)$. That's because, if a new vertex $i^{*}\in {{\cal Z}}_s^c$ is included in this larger cycle, then the outgoing edge from $i^{*}$ in this cycle has a weight with factor $s_{i^{*}}$,  leading to a monomial containing factor $s_{i^{*}}\prod\nolimits_{j = 1}^{{k_{\min }}} {{s_{{{i}_j}}}}$. Next, in the subgraph of ${\mathbb D}^\mathcal{J}$ induced by $\{i_1,...,i_{k_{\min}}\}$,  denoted by ${\mathbb D}_{\bar C_{\min}}^\mathcal{J}$, suppose that there exists another $\lambda$-cycle $\bar C'{\rm{ = }}\{ {i'_1} \to  \cdots  \to {i'_{{k_{\min }}}} \to {i'_1}\}$ which is distinct from $\bar C_{\min}$. Then, according to Lemma \ref{lemma_proof_mid_2}, for every edge $e$ in $E(\bar C_{\min})\bigcup E(\bar C')$, there exists a cycle in ${\mathbb D}_{\bar C_{\min}}^\mathcal{J}$ with length not exceeding $k_{\min}-1$ that contains $e$. Consequently, a $\lambda$-cycle with length less than $k_{\min}$ emerges in ${\mathbb D}^\mathcal{J}$, which is contrary to the shortest length assumption of $\bar C_{\min}$. Therefore, among all the terms ${\small \prod\nolimits_{e \in
{C_i}} {w(e)}}$ which has the form $f(\lambda)\prod\nolimits_{j = 1}^{{k_{\min }}} {{s_{{{i}_j}}}}$ where $f(\lambda)\in F(\lambda)$, there is only one term that corresponds to a matching containing a $\lambda$-edge, while all the rest correspond to matchings that consist of only constant-edges (possibly including self-loops with weight $-1$). These terms have the form of $a\prod\nolimits_{j = 1}^{{k_{\min }}} {{s_{{{i}_j}}}}$ where $a\in \mathbb{R}$ is constant. Hence, the coefficient of  monomial $\prod\nolimits_{j = 1}^{{k_{\min }}} {{s_{{{i}_j}}}}$ in ${\rm{det}} ({G^{\mathcal{J}}_{zv}}{(\lambda )}{P_{\mathcal{J}}} - I)$ has the form \begin{equation}\label{equation_cofficient}  \tilde a+\prod\nolimits_{j = 1}^{{k_{\min }}} {\bar w((i_j,i_{j+1}))},\end{equation} where the assocaited sign is ignored, and $\tilde a\in \mathbb{R}$ is the sum of all possible aforementioned constants $a$. From the definition of $G_{zv}(\lambda)$, $\bar w(e)$ for a $\lambda$-edge $e$ has the form of $\alpha + f(\lambda)$ with $f(\lambda)$ being a nonzero strictly proper fraction, $\alpha \in \mathbb{R}$. Consequently, $\prod\nolimits_{j = 1}^{{k_{\min }}} {\bar w(({i_j},{i_{j + 1}}))}$ always has the form of $\hat \alpha + \hat f(\lambda)$, where $\hat f(\lambda)$ is a nonzero  strictly proper fraction and $\hat \alpha \in \mathbb{R}$, which further leads to (\ref{equation_cofficient}) equaling $\tilde \alpha + \hat \alpha+\hat f(\lambda)$. As a result, there is a term $\hat f(\lambda)\prod\nolimits_{j = 1}^{{k_{\min }}} {{s_{{{i}_j}}}}$ in
${\rm det} (G^\mathcal{J}_{zv}(\lambda)P_{\mathcal{J}}-I)$ depending on both $S_{\cal J}$ and $\lambda$ that cannot be cancelled out by other terms.
It further leads to that $G^\mathcal{J}_{zv}(\lambda)P_{\mathcal{J}}-I$ has at least one zero depending on $S_{\cal J}$.{\footnote{{This can be shown as follows. Let $d_{*}(\lambda)=\det (\lambda I-A_{xx})$, $\hat f(\lambda)=\frac{n_1(\lambda)}{d_1(\lambda)}$, where $n_1(\lambda)$ and $d_1(\lambda)$ are coprime and their degrees $\deg n_1(\lambda)< \deg d_1(\lambda)$. We have proven that $\det (G^\mathcal{J}_{zv}(\lambda)P_{\mathcal{J}}-I)$ has a nonvanished component $(\tilde \alpha+\hat \alpha+\hat f(\lambda))\prod\nolimits_{j = 1}^{{k_{\min }}}{{s_{{{i}_j}}}}+(-1)^{n_s} =\frac{[(\tilde \alpha+\hat \alpha)d_{*}(\lambda)+\frac{n_1(\lambda)d_{*}(\lambda)}{d_1(\lambda)}]\prod\nolimits_{j=1}^{k_{\min}}{s_{i_j}}+(-1)^{n_s}{d_{*}(\lambda)}}{d_{*}(\lambda)}$. Hence, $d_{*}(\lambda)\det (G^\mathcal{J}_{zv}(\lambda)P_{\mathcal{J}}-I)$ has $n$ zeros, where $n$ is the dimension of $A_{xx}$. Suppose that all zeros of $G^\mathcal{J}_{zv}(\lambda)P_{\mathcal{J}}-I$ are independent of $S_{\cal J}$. Then, the zero set of  $d_{*}(\lambda)\det (G^\mathcal{J}_{zv}(\lambda)P_{\mathcal{J}}-I)$ must be $\sigma(A_{xx})$ by setting $s_i|_{i=1}^{n_s}=0$. This further requires that $(\tilde \alpha+\hat \alpha)d_{*}(\lambda)+\frac{n_1(\lambda)d_{*}(\lambda)}{d_1(\lambda)}$ has a zero set $\sigma(A_{xx})$. However, this is impossible, as $\deg \frac{n_1(\lambda)d_{*}(\lambda)}{d_1(\lambda)}< \deg d_{*}(\lambda)$. This contradition shows that there exists at least one zero of $G^\mathcal{J}_{zv}(\lambda)P_{\mathcal{J}}-I$ not belonging to $\sigma(A_{xx})$ for almost every value of $S_{\cal J}$. That is,  $G^\mathcal{J}_{zv}(\lambda)P_{\mathcal{J}}-I$ has at least one zero depending on $S_{\cal J}$.}}}    $\hfill\blacksquare$

{{\emph{Proof of Lemma \ref{lemma_proof_mid_3}:}}} Let $R_1=\{1,...,k\}$, $C_1=\{1,...,k\}$, and $C_2=\{k+1,...,k+q\}$.
 It is obvious that ${G_{zv}}(\lambda ) P  - I$ is of FRGR, as $\det({G_{zv}}(\lambda ) P  - I)$ contains a term $(-1)^k$ which cannot be zeroed out. Hence, for any $j\in C_2$, $[{{G_{zv}}(\lambda ) P  - I},{{G_{zu}}(\lambda )}]_{R_1,C_1\cup C_2\backslash \{j\}}$ is of FRGR.

It remains to prove the case where $j\in C_1$. For this purpose, define the weighted multigraph $\mathbb{D}$ associated with $[{{G_{zv}}(\lambda ) P  - I}\ {{G_{zu}}(\lambda )}]$ in the same way stated in the proof of Lemma \ref{lemma_proof_2}. Moreover, denote the vertex set of $\mathbb{D}$ by ${\cal Z}\cup {\cal U}$ with ${\cal Z}=\{z_1,...,z_k\}$, ${\cal U}=\{u_1,...,u_q\}$. Then, it can be seen that for any $j\in C_1$,  $[{{G_{zv}}(\lambda ) P  - I}\ {{G_{zu}}(\lambda )}]_{R_1,C_1\cup C_2\backslash \{j\}}$ corresponds to the digraph obtained by deleting all the outgoing edges from vertex $z_j$ in $\mathbb{D}$.  According to the reachability of  $z_j$, there is a $u_i\in {\cal U}$ such that a path from $u_i$ to $z_j$ exists in $\mathbb{D}$. Denote $\mathbb{D}^{(i)}$ the subgraph of $\mathbb{D}$ induced by vertices $\{u_i,z_1,...,z_k\}$. Let the path ${\cal P}{\rm{\triangleq \{ }}{u_i} \to {z_{{i_1}}}\to\cdots\to{z_{{i_{d - 1}}}}\to{z_j}{\rm{\} }}$ be the shortest path from $u_i$ to $z_j$ in $\mathbb{D}^{(i)}$, whose length is $d$, where $\{i_1,...,i_{d-1}\}\subseteq [k]$. %(if $d=1$, $\{i_1,...,i_{d-1}\}=\emptyset$)
  Let $\mathbb{D}^{(i)}_d$ be the subgraph of $\mathbb{D}^{(i)}$ induced by vertices $\{u_i, z_{i_1},...,z_{{i_{d-1}}}, z_j\}$. Similarly to the proof of Lemma \ref{lemma_proof_2}, it can be seen that the path $\cal P$ and all self-loops with weights $-1$ associated with each vertex of ${\cal Z}\backslash\{z_{i_1},...,z_{{i_{d-1}}}, z_j\}$ constitute a matching with size $k$, denoted by $C_{\cal{P}}$. This matching corresponds to a term with the monomial ${s_{{i_1}}}{s_{{i_2}}}...{s_{{i_{d - 1}}}}$ in $\det([{{G_{zv}}(\lambda ) P  - I}\ {{G_{zu}}(\lambda )}]_{R_1,C_1\backslash \{j\}\cup \{i+k\}})$. It is declared that among all matchings with size $k$ in $\mathbb{D}^{(i)}$, except $C_{\cal{P}}$,  none corresponds to a term with the monomial $s_{i_1}s_{i_2}...s_{i_{d-1}}$. We can demonstrate this by contradiction. If such matching exists, the only possible case is that, there is a path from $u_i$ to $z_j$ with length $d$ in $\mathbb{D}_d^{(i)}$ other than $\cal{P}$. However, by Lemma \ref{lemma_proof_mid_2}, this means that a path from $u_i$ to $z_j$ exists with length no more than $d-1$, causing a contradiction to the shortest path assumption. To show this, one just needs to virtually add an edge $(z_j,u_i)$ to $\mathbb{D}_d^{(i)}$. Then a path from $u_i$ to $z_j$ always corresponds to a cycle in Lemma \ref{lemma_proof_mid_2}. Hence, the monomial ${s_{{i_1}}}{s_{{i_2}}}...{s_{{i_{d - 1}}}}$ cannot be canceled in  $\det ([{{G_{zv}}(\lambda ) P  - I}\ {{G_{zu}}(\lambda )}]_{R_1,C_1\backslash \{j\}\cup \{i+k\}})$, making $[{{G_{zv}}(\lambda ) P  - I}\ {{G_{zu}}(\lambda )}]_{R_1,C_1\cup C_2\backslash \{j\}}$ FRGR. $\hfill\blacksquare$

{\emph{Proof of Lemma \ref{lemma_proof_mid_4}:}} The proof follows similar ideas to those of \citep[Lemma 5.3]{Murota_SIAM}.  Suppose that there exists a zero of $[{{G_{zv}}(\lambda ) P\! -\! I}\ {{G_{zu}}(\lambda )}]$ depending on $S$, denoted by $\lambda^*$. Then, inspired by \citep[Lemma 5.3]{Murota_SIAM}, it suffices to see $\lambda^*$ as a transcendental element over $\mathbb{R}$, since $\lambda^*$ depends on the algebraically independent elements $s_1,...,s_k$. Choosing the square submatrix ${{G_{zv}}(\lambda ) P  - I}$ of $[{{G_{zv}}(\lambda ) P  - I}\ {{G_{zu}}(\lambda )}]$, it holds that $\det ({{G_{zv}}(\lambda^* ) P  - I})=0$ since $\lambda^*$ is the zero, which means that $\{\lambda^*\}\cup \{s_1,...,s_k\}$ is algebraically dependent over $\mathbb{R}$. According to the property of algebraic independence \citep[Lemma 2.1]{Murota_SIAM}, there exists some $s_{i^*}\in \{s_1,...,s_k\}$ such that $\{\lambda^*\}\cup \{s_1,...,s_k\} \backslash \{s_{i^*}\}$ is algebraically independent over $\mathbb{R}$. On the other hand, by Lemma \ref{lemma_proof_mid_3}, for any $i^*\in [k]$, deleting the $i^*$-th column from $[{{G_{zv}}(\lambda ) P \! - \!I}\ {{G_{zu}}(\lambda )}]$, the obtained matrix is of FRGR. This means that, at least one nonsingular $k\times k$ submatrix of $[{{G_{zv}}(\lambda ) P  \!-\! I}\ {{G_{zu}}(\lambda )}]$ does not contain $s_{i^*}$. Denote the determinant of such submatrix by $D(S\backslash \{s_{i^*}\},\lambda^*)$, then $D(S\backslash{\{s_{i^*}\}},\lambda^*)=0$ holds. It follows that $\{\lambda^*\}\cup \{s_1,...,s_k\} \backslash \{s_{i^*}\}$ is algebraically dependent over $\mathbb{R}$, causing a contradiction. $\hfill\blacksquare$

{{\emph{Proof of Proposition \ref{proposition_3}:}}} Define a digraph associated with $[P G_{zv}(\lambda)\ P G_{zu}(\lambda)]$ as $\hat {\mathbb{T}}_{{\bf{{\bf \Sigma}}}}=({\cal{V}} \cup {\cal{U}}, {\cal{E}}_{{\cal V}{\cal V}}\cup {\cal{E}_{{\cal U}\cal{V}}})$, where ${{\cal{U}}}\triangleq \{u_1,...,u_{M_u}\}$, ${{\cal{V}}}\triangleq \{v_1,...,v_{M_v}\}$, ${\cal{E}}_{{\cal V}{\cal V}}=\{(v_i,v_j): [PG_{zv}(\lambda)]_{ji}\ne 0\}$ and ${\cal{E}}_{{\cal U}{\cal V}}=\{(u_i,v_j): [P G_{zu}(\lambda)]_{ji}\ne 0\}$. {{Noting that all nonzero entries of $P$ are independent, it turns out that $[P G_{z\star}(\lambda)]_{ji}\ne 0$, if and only if there exists $l\in[M_z]$ such that $P_{jl}\ne 0$ and $[G_{z\star}(\lambda)]_{li}\ne 0$, where $\star = v$ or $u$. Relabel vertices in $\mathbb{T}_{{\bf{\Sigma}}}$ as $\bigcup\nolimits_{i=1}^{N}{{\cal{U}}_i}\triangleq \{u_1,...,u_{M_u}\}$, $\bigcup\nolimits_{i=1}^N{{\cal{Z}}_i}\triangleq \{z_1,...,z_{M_z}\}$ and $\bigcup\nolimits_{i=1}^N{{\cal{V}}_i}\triangleq \{v_1,...,v_{M_v}\}$ according to their correspondences in $\hat {\mathbb{T}}_{{\bf{{\bf \Sigma}}}}$. Then, $\mathbb{T}_{{\bf{\Sigma}}}$ and $\hat {\mathbb{T}}_{{\bf{{\bf \Sigma}}}}$ are related in the following way:  $(v_i,v_j)\in {\cal{E}}_{\cal{V}\cal{V}}$ (resp. $(u_i,v_j)\in {\cal{E}_{\cal{U}\cal{V}}}$) in $\hat {\mathbb{T}}_{{\bf{{\bf \Sigma}}}}$, if and only if there is a triple $(v_i,z_{i^*},v_j)$ for some $i^*$, such that $(v_i,z_{i^*}), (z_{i^*},v_j)\in \cal{E}_{\bf{\Sigma}}$ (resp., a triple $(u_i, z_{i^*}, v_j)$ with $(u_i,z_{i^*}), (z_{i^*},v_j)\in \cal{E}_{\bf{\Sigma}}$) in $\mathbb{T}_{{\bf{{\bf{\Sigma}}}}}$. It can further be validated that, there exists an input-unreachable $\lambda$-cycle in $\mathbb{T}_{{\bf{{\bf \Sigma}}}}$, if and only if there exists at least one input-unreachable $\lambda$-cycle in $\hat {\mathbb{T}}_{{\bf{{\bf \Sigma}}}}$.

  Define $H(\lambda)$ and $G(\lambda)$ as $H(\lambda)\triangleq V{G_{zu}(\lambda)}$, $G(\lambda)\triangleq V{G_{zv}(\lambda)}$ respectively. Noting that $P=UP_dV$ and $P_d \triangleq {\rm{{\bf{{\bf {\bf diag}}}}}}\{s_i|_{i=1}^{k}\}$, it can be seen that $[PG_{zv}(\lambda)]_{ij}\ne 0$, i.e., $[UP_dVG_{zv}(\lambda)]_{ij}\ne 0$ (resp. $[PG_{zv}(\lambda)]_{ij}$ depends on $\lambda$), if and only if there exists an $l\in[k]$, such that $U_{il}\ne0$ and $[P_dG(\lambda)]_{lj}\ne 0$, where the latter is equivalent to $[G(\lambda)]_{lj}\ne 0$ (resp. $[G(\lambda)]_{lj}$ depends on $\lambda$). Similar analysis is valid for $PG_{zu}(\lambda)$. Hence, it turns out that $[PG_{zv}(\lambda)\  PG_{zu}(\lambda)]$ has the same sparsity pattern as $\left[U\!*\!G(\lambda)\ U\!*\!H(\lambda)\right]$, where the operation  ``$\ast$'' is defined as the same as matrix multiplication except that the involved additions and multiplications between two scalar elements are logic operations OR and AND between binaries respectively (here binarity refers to whether an entry is zero or nonzero). Furthermore, since each column of $U$ has only one nonzero entry from its definition, no cancellation occurs between two addends in obtaining $G(\lambda)U$. This means that $[VG_{zv}(\lambda)U\ VG_{zu}(\lambda)]$ has the same sparsity pattern as $\left[G(\lambda)\!*\!U\ H(\lambda)\right]$. With these observations, we will show that $\hat {\mathbb{T}}_{{\bf{{\bf \Sigma}}}}$ and $\mathbb{L}_{{\bf{{\bf \Sigma}}}}$ have the same properties w.r.t. the existence of $\lambda$-cycles and their input-reachabilities by leveraging the correspondence in sparsity pattern between their associated  matrices $\left[U\!*\!G(\lambda)\ U\!*\!H(\lambda)\right]$ and $\left[G(\lambda)\!*\!U\ H(\lambda)\right]$.}}

To this end, denote the set of vertices of $\mathbb{L}_{{\bf{{\bf \Sigma}}}}$ by $\cal{W}\cup \cal{U}$ with ${\cal{W}}=\{w_1,...,w_k\}$ and ${\cal{U}}=\{u_1,...,u_{M_u}\}$.  Suppose that there is a $\lambda$-cycle in $\mathbb{L}_{{\bf{{\bf \Sigma}}}}$, denoted by ${{\cal{C}}_1} \triangleq  \{{w_{{i_1}}}\!\to\! {w_{{i_2}}} \!\to\!  \cdots  \!\to\! {w_{{i_s}}} \!\to\! {w_{{i_1}}}\}$. Assume that $(w_{i_1},w_{i_2})$ is a $\lambda$-edge without loss of generality. Moreover, suppose that there is a path from $u_{\bar i_0}\in {\cal{U}}$ to $w_{\bar i_r}\in \{w_{i_1},...,w_{i_s}\}$, and denote such path by $\{{u_{\bar i_0}} \!\to\! {w_{\bar i_1}} \!\to\!  \cdots  \!\to\! {w_{\bar i_r}}\}$.  This means that the $({\bar i_1,\bar i_0})$-th entry of $H(\lambda)$, and the $(\bar i_2,\bar i_1)$-th,...,$(\bar i_r,\bar i_{r\!-\!1})$-th,...,$({i_2},{i_1})$-th,...,$({i_s},{i_{s\!-\!1}})$-th and $({i_1},{i_s})$-th entries of $G(\lambda)*U$ are nonzeros. {{ Note that $[G(\lambda)*U]_{ij}\ne 0$, if and only if there exists an $l\in[M_v]$ such that $[G(\lambda)]_{il}\ne 0$ and $U_{lj}\ne 0$. Hence, there exists a sequence of integers $\bar k_1,...,\bar k_{r\!-\!1},{k_1},...,{k_s} \in [M_v]$ (possibly with repeated values), such that $[G(\lambda)]_{\bar i_{j\!+\!1}, \bar k_{j}}\ne 0$ and $U_{\bar k_{j}, \bar i_{j}}\ne 0$ for $j=1,...,r-1$, $[G(\lambda)]_{i_{j\!+\!1}, k_{j}}\ne 0$ and $U_{k_{j}, i_j}\ne 0$ for $j=1,...,s$, where $i_{s+1}$ is defined to be $i_1$, and meanwhile $[G(\lambda)]_{i_2,k_1}$ depends on $\lambda$. Observe that whenever $[G(\lambda)]_{il}\ne 0$ and $U_{ji}\ne 0$, $[U\!*\!G(\lambda)]_{jl}\ne 0$. The above indicates that the $(\bar k_2,\bar k_1)$-th,...,$(k_{\bar s},\bar k_{r\!-\!1})$-th,$(k_1,k_s)$-th,$(k_2,k_1)$,...,$(k_{s},k_{s-1})$-th entries of $U\!*\!G(\lambda)$ are nonzeros, where $\bar s\in [s]$, and $[U\!*H(\lambda)]_{\bar k_1, \bar i_0}\ne 0$. Then, there exists a sequence of edges $(u_{\bar i_0},v_{\bar k_1}),...,(v_{\bar k_{r\!-\!1}},v_{k_{\bar s}})$,...,$(v_{k_{s}},v_{k_{1}})$, $(v_{k_1},v_{k_2})$,...,$(v_{k_{s\!-\!1}},v_{k_{s}})$ in $\hat {\mathbb{T}}_{{\bf \Sigma}}$ (possibly with repeated edges) with $(v_{k_1},v_{k_2})$ being a $\lambda$-edge.}} A $\lambda$-cycle containing edge $(v_{k_1},v_{k_2})$ can always be found from these edges, and this $\lambda$-cycle is input-reachable.

Since every step of the above analysis is invertible, the direction that a (input-reachable) $\lambda$-cycle in $\hat {\mathbb{T}}_{{\bf{{\bf \Sigma}}}}$ indicates the existence of a (input-reachable) $\lambda$-cycle in $\mathbb{L}_{\bf{{\bf \Sigma}}}$ follows a similar way. This further leads to Proposition \ref{proposition_3}.  $\hfill\blacksquare$

{{\emph{Proof of Lemma \ref{Feasibility1}:}}}
The necessity of (i) can be directly derived from Theorem 1 of \cite{Y_Zhang_2016}, \cite{zhou_2018}. The necessity of (ii) can be validated by contradiction: if $M_z< \max\nolimits_{1\le i \le m}{M_{ri}}$, $\rho({\cal M}(Q_1)\cap {\cal M}(Q_2)) \le M_z < \max\nolimits_{1\le i \le m}{M_{ri}}$, which never satisfies the condition of Proposition \ref{proposition 4}. The necessity of (iii) is a direct derivation of Condition (1) of Corollary \ref{corollary_add} and the definition of ${\mathbb{T}}_i$.

To show the sufficiency, let $\Phi_{{\rm full}}$ be the SCM with all of its entries being free parameters. Notice that if (i) of Lemma \ref{Feasibility1} is satisfied,
$[Y^{(j)}_{i}\ Z^{(j)}_{i}]$ is of FRR whenever $m_{rij}>0$ by using Lemma \ref{lemma3} inversely, which means that $Q_{2i}$ is of FRR.  Suppose $M_z\ge M_{ri}$ for all $i\in [m]$. Then, any $M_{ri}$ columns of $Q_1 \triangleq [\Phi_{{\rm full}}^{\intercal}\ I_{M_z}]$ are linearly independent. As $Q_{2i}$ is of FRR, there exists $M_{ri}$ columns of $Q_{2i}$ which are linearly independent. As a result, $\rho({\cal M}({Q_1}) \cap {\cal M}({Q_{2i}})) = {M_{ri}}$ per $i\in [m]$. If there is some $i\in [N]$ with $G^{(i)}_{zu}(\lambda)\ne 0$, it can be validated that $\Phi_{{\rm full}}$ is sufficient to make all $\lambda$-edges in the associated $n$-ACG ${\mathbb{T}}_{\Sigma}$ input-reachable.  Therefore, $\Phi_{{\rm full}}$ satisfies both conditions of Corollary \ref{corollary_add}, which is a feasible solution to Problem 2. $\hfill\blacksquare$

{{\emph{Proofs of Proposition \ref{NP_hard1} and Lemma \ref{NP_hard2} :}}} We put both proofs together as they follow the same argument. The sketch is to find an instance of Problem 2 (resp. Problem \ref{prob3}) that is equivalent to the NP-hard {{minimal controllability problem}} in \cite{A_Olsehvsky_2014}. The latter problem is to determine the minimal number of states that need to be affected by an input to ensure controllability for a given state transition matrix \cite{A_Olsehvsky_2014}. Consider an NDS ${\bf{{\bf \Sigma}}}$ with two subsystems ${\bf{{\bf \Sigma}}}_1$ and ${\bf{{\bf \Sigma}}}_2$. The parameters are as follows:  for a given $n\in \mathbb{N}$, let $A^{(1)}_{xx}\in \mathbb{R}^{n\times n}$ be a matrix with no repeated eigenvalues whose associated minimal controllability problem is NP-hard, whose construction can be referred to Theorem 1 of \cite{A_Olsehvsky_2014}, $A^{(1)}_{xv}=I_n$, $A^{(1)}_{zv}=0_{n\times n}$, $A^{(1)}_{zx}=I_n$, $B^{(1)}_{zu}=0_{n\times 1}$, $B^{(1)}_{xu}=0_{n\times 1}$;  $A^{(2)}_{xx}=n+1$, $A^{(2)}_{xv}=1$, $A^{(2)}_{zv}={{0}}_{1\times 1}$, $A^{(2)}_{zx}=1$, $B^{(2)}_{zu}=0_{1\times 1}$, $B^{(2)}_{xu}=1$, and the SCM to be determined $\Phi\in\{0,*\}^{(n+1)\times (n+1)}$. Then, it can be validated that, Problem 2 and Problem \ref{prob3} are both  equivalent to the minimal controllability problem associated with $A^{(1)}_{xx}$. $\hfill\blacksquare$

{{\emph{Proof of Proposition \ref{Submodular_proposition}:}}} Define functions $g_i(J)\triangleq {\rm rank}([Y_{iJ}\  Z_i])$ per $i\in [m]$. As $g_i(J)$ is a rank function on the subset of column vectors of matrix $[Y_{i}\   Z_i]$, $g_i(J)$ is submodular and nondecreasing on $J\subseteq [M_v]$ \cite{Murota_Book}.
Hence, $g(J)=\sum \nolimits_{i=1}^{m} g_i(J)$ is submodular and nondecreasing.
Denote the optimal solution to Problem \ref{prob4} by $J^*$. From \cite{Submodular}, it follows that{\small{
\begin{equation} \label{Submodular_mid} \frac{{|{J_{\rm grd}}|}}{{\left| {{J^*}} \right|}} \le 1 + \log \frac{{\sum\nolimits_{i = 1}^{m} {{M_{ri}} - \sum\nolimits_{i = 1}^{m} {{\rm rank}({Z_i})} } }}{{\sum\nolimits_{i = 1}^{m} {{M_{ri}} - g({J_{T - 1}})} }}.\end{equation}}}
Let $R_1^{*}\subseteq [M_v]$ be the set of indices of nonzero rows of $\Phi^{*}$. Then, it is obvious that $||\Phi^* ||_{0}\ge |R_1^{*}|$. Let ${\cal{B}}\triangleq \{M_{v}+1,..., M_{v}+M_{z}\}$. Since $\Phi^{*}$ satisfies $f(\Phi^*)=\sum \nolimits_{i = 1}^{m} {M_{ri}}$, it means that ${\rho ({\cal{M}}({Q_1}(\Phi^* )) \cap {\cal{M}}(Q_{2i}))}=M_{ri}$ per $i\in [m]$. That is, there is a collection $\{J_1,..., J_{m}\}$ of subsets of $[M_v+M_z]$, such that $|J_i|=M_{ri}$, ${\rm rank}(Q_{2i,J_i})=M_{ri}$ and ${\rm rank}(Q_1(\Phi^*)_{J_i})=M_{ri}$.  It further leads to that $J_i\backslash {\cal B} \subseteq R_1^{*}$, and ${\rho ({\cal M}(\Omega (J_i\backslash {\cal B})) \cap {\cal M}({Q_{2i}}))}=M_{ri}$. Moreover, define $R_2^* \triangleq  \bigcup\nolimits_{i = 1}^{m} {({J_i}\backslash {\cal B})}$. Then, the above relations imply that $R_2^*\subseteq R_1^*$ and  $g(R_2^*)=\sum\nolimits_{i=1}^{m}{M_{ri}}$, i.e., $R_2^*$ is feasible for Problem \ref{prob4}. Consequently, we have $|J^*|\le |R^{*}_2|\le |R^{*}_1|\le ||\Phi^* ||_{0}$.

By the definition of $Z^{(j)}_{i}$, per $i\in [m]$, $j\in [N]$, it holds
$m_{rij}-{\rm rank}(Z^{(j)}_{i})=m_{x{j}}- {\rm rank}([\lambda_i I- A_{xx}^{(j)}\ B_{xu}^{(j)}])\le m_{xj}-{\rm rank}(\lambda_i I-A^{(j)}_{xx})$,  where the last term is the geometric multiplicity of $\lambda_i$ for $A^{(j)}_{xx}$.
 Combined with (\ref{Submodular_mid}), summation of the above relations over $[m]$ and  $[N]$ leads to the inequalities of Proposition \ref{Submodular_proposition}.  $\hfill\blacksquare$

{{\emph{Proof of Theorem \ref{theorem_last}:}}}
Let $\Phi^*$ denotes the optimal solution to Problem \ref{prob3},  $\Phi'$ the returned solution by Algorithm \ref{alg3}, and $\Phi_{\rm opt}$ the optimal solution to Problem \ref{prob2}.

We first show the feasibility. For the obtained $\mathbb{G}_{\rm col}$ in the last iteration of Step 2 of Algorithm \ref{alg3}, we say a vertex is $k$-colored, if it is assigned with $k$ different colors, $1\le k \le M_{\rm rmax}$.  Let every vertex of $\mathbb{G}(J_1,...,J_{m})$ has the same colors as the corresponding vertex of $\mathbb{G}_{\rm col}$. Considering the subgraph of $\mathbb{G}(J_1,...,J_{m})$ induced by $J_i$ per $i\in [m]$, denoted by $\mathbb{G}_{J_i}$, each of its vertex is either $1$-colored or $k^{+}$-colored for some $k^+>1$. Recall that a vertex is $k^{+}$-colored only if it has at least $M_z$ differently $1$-colored neighbors. Hence, all vertices in ${J_i}\cap {\cal B}$ are $1$-colored. According to (\ref{color_k}) in Algorithm \ref{alg3}, all those $k^+$ satisfies $k^+ \ge M_{ri}$.   As a result, there always exists one combination of $M_{ri}$ colors in $\mathbb{G}_{J_i}$, denoted by $J'_i\subseteq [M_z]$ with $|J'_i|=M_{ri}$, each color chosen from each vertex separately, such that $\mathbb{G}_{J_i}$ is colored with the property that no two adjacent vertices share the same chosen color (for example, first choose the unique colors from those $1$-colored vertices, then combinatorially choose one color from each of the $k^+$-colored ($k^+\ge M_{ri}$) vertices). By the relation between coloring and rank, it follows that $Q_1(\Phi')_{J'_i,J_i}$ has full generic rank. As ${\rm rank}(Q_{{2i},J_i})= M_{ri}$, we have ${\rho ({\cal{M}}({Q_1}(\Phi')) \cap {\cal{M}}(Q_{2i}))}=M_{ri}$.  This means that $\Phi'$ is feasible for Problem \ref{prob3}. The feasibility of the two-stage algorithm then follows from Propositions \ref{proposition_3} and \ref{proposition 4}.

We then prove these approximation bounds. From the above analysis, it is clear that every vertex of $\mathbb{G}_{\rm col}$ is colored with no more than $M_{\rm rmax}$ colors.  Combining Proposition \ref{Submodular_proposition},  we have $||\Phi'||_0/||\Phi^*||_0 \le M_{\rm rmax}({\rm log}(M_{\rm def})+1)$. Moreover, if there does not exist a $k^{+}$-colored vertex in $\mathbb{G}_{\rm col}$ for any $k^{+}>1$, then obviously $||\Phi'||_0/||\Phi^*||_0 \le {\rm log}(M_{\rm def})+1$.  For the overall
topology design procedure,  by Corollary \ref{corollary_add}, it can be seen that subsystem links in an optimal solution $\Phi_{\rm opt}$ to  {{Problem \ref{prob2}}} can be divided into two subsets. One subset functions as eliminating the FUMs, which is supposed to contain at least $p^{*}_{\rm ef}$ subsystem links, and the other functions as eliminating the input-unreachable $\lambda$-edges, which is supposed to contain at least $p^{*}_{\rm eu}$ subsystem links (these two subsets may overlap). It turns out that both $p^{*}_{\rm ef}$ and $p^{*}_{\rm eu}$ do not increase with the addition of any set of subsystem links to the corresponding NDS.  Hence,  we have
$$p^{*}_{\rm ef}+p^{*}_{\rm eu}\le 2||\Phi_{\rm opt}||_0, ||\Phi^*||_0\le p^{*}_{\rm ef}.$$
On the other hand, given a collection of {{disconnected}} subsystems ${\bf \Sigma}_1$,...,${\bf \Sigma}_N$, it holds that
$p_{\rm {\rm {ius}}}\le p^{*}_{\rm eu},$ recalling that $p_{\rm {\rm {\rm {ius}}}}$ is the number of input-unreachable source SCCs that contains an $\lambda$-edge in ${\mathbb{T}}_{{\bf{{\bf \Sigma}}}{\text{s1}}}$.  This is obvious, as the existence of $p_{\rm {\rm {\rm {ius}}}}$ input-unreachable $\lambda$-edges in different source SCCs of ${\mathbb{T}}_{{\bf{{\bf \Sigma}}}{\text{s1}}}$, indicates that at least $p_{\rm {\rm {\rm {ius}}}}$ subsystem links should be added to $\bigcup\nolimits_{i=1}^N {{{\mathbb{T}}}_i}$ to make the corresponding input-unreachable $\lambda$-edges therein input-reachable. Hence, the two-stage algorithm returns a $\Phi$ with sparsity
$p_{\rm {ius}}+||\Phi'||_0\le p^{*}_{\rm eu}+ M_{\rm rmax}({\rm log}(M_{\rm def})+1)p^{*}_{\rm ef} \le M_{\rm rmax}(1+{\rm log}(M_{\rm def}))(p^{*}_{\rm eu}+p^{*}_{\rm ef})
\le 2M_{\rm rmax}({\rm log}(M_{\rm def})+1)||\Phi_{\rm opt}||_0$. Therefore, the topology design procedure overall has an $O(2M_{\rm rmax}{\rm log}(M_{\rm def}))$ approximation.
$\hfill\blacksquare$

\end{appendices}


\begin{thebibliography}{20}
\footnotesize

\bibitem{Anderson_1982}
B. D. O. Anderson and H. Hong, “Structural controllability and
matrix nets,” \emph{Int. J. Contr.}, vol. 35, no. 3, pp. 397-416, 1982.

\bibitem{partitions}
R. Amirreza, M. Ji, M. Mesbahi, and M. Egerstedt, controllability of multi-agent systems from a graph-theoretic perspective, \emph{SIAM J. Contr. Optim.}, vol. 48, no. 1, pp. 162-186, 2009.

\bibitem{computation_complexity_modern_2009}
S. Arora and B. Boaz, \emph{Computational Complexity: a Modern Approach}, Cambridge University Press, 2009.

\bibitem{network of networks}
A. Chapman, M. N. Abdolyyousefi, and M. Mesbahi, Controllability
and observability of network-of networks via
cartesian products,  \emph{IEEE Trans. Automat. Contr.}
vol. 59, no. 10, pp. 2668-2679, 2014.

\bibitem{Clark_input_selection}
A. Clark, B. Alomair, L. Bushnell, et al, “Input selection for performance and controllability of structured linear descriptor systems,”  \emph{SIAM J. Contr. Optim.}, vol. 55, no. 1, pp. 457-485, 2017.

\bibitem{Morse_1976}
J. Corfmat and A. S. Morse, “Structurally controllable and structurally
canonical systems,” \emph{IEEE Trans. Automat. Contr.}, vol. 21, no. 1, pp.
129-131, 1976.

\bibitem{Morse_decentralized}
J. Corfmat and A. S. Morse, Decentralized control of linear multivariable systems, \emph{Automatica}, vol. 12, no. 5,  pp. 479-495, 1976.

\bibitem{Composability}
J. F. Carvalho, S. Pequito,  A. P. Aguiar, et al., Composability and controllability of structural linear time-invariant systems: distributed verification, \emph{Automatica}, vol. 78, pp. 123-134, 2017.

\bibitem{plos one}
N. J. Cowan, E. J. Chastain, D. A. Vilhena, et al., Nodal dynamics, not degree distributions, determine the structural controllability of complex networks, {\emph{PloS One}}, vol. 7, no. 6, e38398, 2012.



\bibitem{generic}
J. M. Dion, C. Commault, and J. Van DerWoude, Generic properties and control of linear structured
systems: a survey, {\emph{Automatica}}, vol. 39, pp. 1125-1144, 2003.




\bibitem{Matrix analysis}
R. A. Horn and C. R. Johnson, \emph{Topics in Matrix Analysis}, Cambridge
University Press, 1991.



\bibitem{Derandom_2004}
V. Kabanets and R. Impagliazzo, Derandomizing polynomial identity tests means proving circuit lower bounds, \emph{Comp. Comple.}, vol. 13, pp. 1-46, 2004.



\bibitem{Kailath_1980}
T. Kailath, \emph{Linear Systems}, Englewood Cliffs, NJ: Prentice Hall, 1980.

\bibitem{S_Pequito_2017}
 S. Kruzick, S. Pequito, S. Kar, et al, Structurally observable distributed networks of agents under cost and robustness constraints,  \emph{IEEE Trans. Sig. Info. Pr. Netw.} , DOI: 10.1109/TSIPN.2017.2681208, 2017.

\bibitem{langbort_2004}
C. Langbort, R. S. Chandra, and R. D'Andrea,  Distributed control design for systems interconnected over an arbitrary graph, \emph{IEEE Trans. Automat. Contr.},  vol. 49, no. 9, pp. 1502-1519, 2004.

\bibitem{Liu_Fei}
F. Liu and A. S. Morse,  Structural controllability of linear time-invariant systems, {\emph{arXiv preprint arXiv}}:1707.08243, 2017.


% \bibitem{Rational_function}
% K. S. Lu and J. N. Wei, Rational function matrices and structural controllability and observability, \emph{IEE Proc. D}, vol. 138. no. 4, 1991.

\bibitem{Lin_1974}
C. T. Lin, Structural controllability, \emph{IEEE Trans. Automat. Contr.}, vol. 19, no. 3, pp. 201-208, 1974.

\bibitem{nature}
Y. Y. Liu, J. J. Slotine, and A. L. Barabasi, Controllability of complex networks, \emph{Nature}, vol. 473, no.7346, pp.167-
173, 2011.


\bibitem{Rational_function}
K. S. Lu and J. N. Wei, Rational function matrices and structural controllability and observability, \emph{IEE Proc. D}, vol. 138. no. 4, pp. 388-394, 1991.

\bibitem{Murota_SIAM}
K. Murota,  “Refined study on structural controllability of descriptor systems by means of matroids,” \emph{SIAM J. Contr. Optim.},  vol. 25, no. 4,  pp. 967-989, 1987.

\bibitem{Murota_Book}
K. Murota, \emph{Matrices and Matroids for Systems Analysis}, Springer Science Business Media, 2009.

\bibitem{strong_controllability}
H. Mayeda and T. Yamada, Strong structural controllability, \emph{ SIAM J. Contr. Optim.}, vol. 17, no. 1, pp. 123-138, 1979.

\bibitem{A_Olsehvsky_2014}
A. Olshevsky, Minimal controllability problems, {\emph{IEEE Trans. Contr. Contr. Netw. Syst.} }, vol.1, no.3, pp. 249-258, 2014.

\bibitem{Modern_Control_Ogata}
K. Ogata and Y. Yang,  \emph{Modern Control Engineering}, India: Prentice Hall, 2002.

\bibitem{controllability metrics}
F. Pasqualetti, S. Zampieri, and F. Bullo, Controllability metrics, limitations and algorithms for complex networks, \emph{IEEE Trans. Contr. Netw. Syst.}, vol.1, no.1, pp. 40-52, 2014.

\bibitem{FBullo2013}
 F. Pasqualetti, F. Dorfler, and F. Bullo, Attack detection and identification in cyber-physical systems, \emph{IEEE Trans. Automat. Contr.}, vol. 58, no. 11, 2013.


\bibitem{A_Shpilka_survey}
A. Shpilka and A. Yehudayoff, Arithmetic circuits: A survey of recent results and open questions, {\emph{Found. Tren. Theor. Comput. S.}}, vol. 5, pp. 207-388, 2010.

\bibitem{Separation_1980}
J. T. Schwartz, Fast probabilistic algorithms for verification of polynomial identities, \emph{J. Assoc. Comput. Mach.}, vol. 27, pp. 701-717, 1980.

\bibitem{Tsopelakos bilinear}
A. Tsopelakos,  M. A. Belabbas, and B. Gharesifard,  Classification of the structurally
controllable zero-patterns for driftless bilinear control systems, \emph{IEEE Trans. Contr. Netw. Syst.}, doi: 10.1109/TCNS.2018.2834822, 2018.


\bibitem{Willems_1986}
J. L. Willems, Structural controllability and observability, \emph{Syst. ${\rm{\& }}$
Control Lett.}, vol. 8, no. 1, pp. 5-12, 1986.

\bibitem{Submodular}
L. A. Wolsey, An analysis of the greedy algorithm for the submodular set covering problem, \emph{Combinatorica},
vol. 2, no. 4, pp. 385-393, 1982.

\bibitem{DB_West_graph}
D. B. West, \emph{Introduction to Graph Theory}, Upper Saddle River: Prentice Hall, 2001.


\bibitem{Fixed Mode}
S. H. Wang and E. J. Davison,  On the stabilization of decentralized control systems,  \emph{IEEE Trans. Automat. Contr.},  vol. 18, no. 5,  pp. 473-478, 1973.



\bibitem{Zhou_robust_book}
K. Zhou, J. C. Doyle, and K. Glover, \emph{Robust and Optimal Control}, New Jersey: Prentice Hall, 1996.


\bibitem{zhou_2015}
T. Zhou, On the controllability and observability of networked dynamic systems, \emph{Automatica}, vol. 52, pp. 63-75, 2015.

\bibitem{zhou_2018}
T. Zhou, Minimal inputs outputs for subsystems in a networked system, \emph{Automatica}, vol. 94, pp. 161-169, 2018.

\bibitem{Y_Zhang_2016}
Y.  Zhang and T. Zhou, Controllability analysis for a networked dynamic
system with autonomous subsystems, \emph{IEEE Trans. Automat. Contr.}, vol. 62, no. 7, pp. 3048-3415, 2017.

\bibitem{Y_Zhang_2017}
Y. Zhang and T. Zhou, On the edge insertion/deletion and
controllability distance of linear structural systems, in \emph{Proc. IEEE Conf. Dec. Control}, pp. 2300-2305, 2017.

\bibitem{Y_Zhang_2018}
Y. Zhang and T. Zhou, Input matrix construction and approximation using a graphic approach, {\emph{Int. J. Contr.}}, doi: 10.1080/00207179.2018.1519600, 2018.
\end{thebibliography}
\end{document}